# OPTIMAL RANK-BASED TESTS FOR HOMOGENEITY OF SCATTER

By Marc Hallin[1] and Davy Paindaveine[2]

*Université Libre de Bruxelles*

We propose a class of locally and asymptotically optimal tests, based on multivariate ranks and signs for the homogeneity of scatter matrices in $m$ elliptical populations. Contrary to the existing parametric procedures, these tests remain valid without any moment assumptions, and thus are perfectly robust against heavy-tailed distributions (*validity robustness*). Nevertheless, they reach semiparametric efficiency bounds at correctly specified elliptical densities and maintain high powers under all (*efficiency robustness*). In particular, their normal-score version outperforms traditional Gaussian likelihood ratio tests and their pseudo-Gaussian robustifications under a very broad range of non-Gaussian densities including, for instance, all multivariate Student and power-exponential distributions.

## 1. Introduction.

1.1. *Homogeneity of variances and covariance matrices.* The assumption of variance homogeneity is central to the theory and practice of univariate $m$-sample inference, playing a major role in such models as $m$-sample location (ANOVA) or $m$-sample regression (ANOCOVA). The problem of testing the null hypothesis of variance homogeneity, therefore, is of fundamental importance, and for more than half a century, has been a subject of continued interest in the statistical literature. The standard procedure, described in most textbooks, is Bartlett's modified (Gaussian) likelihood ratio test (see [2]). This test, however, is well known to be highly nonrobust against violations of Gaussian assumptions, a fact that gave rise to a large number

Received May 2007; revised May 2007.
[1]Supported by the Bendheim Center.
[2]Supported by a Crédit d'Impulsion of the Belgian Fonds National de la Recherche Scientifique.

of "robustified" versions of the likelihood ratio procedure ([3, 5, 6, 21], to quote only a few). Soon, it was noticed that these "robustifications," if reasonably resistant to nonnormality, unfortunately were lacking power: in the convenient terminology of Heritier and Ronchetti [22], they enjoy *validity robustness* but not *efficiency robustness*.

In an extensive simulation study, Conover, Johnson and Johnson [7] have investigated the validity robustness (against nonnormal densities) and efficiency robustness properties of 56 distinct tests, including several (signed-) rank-based ones. Their conclusion is that only three of them survive the examination, and that two of the three survivors are normal-score signed-rank tests (adapted from [10]).

In view of its applications in MANOVA, MANOCOVA, discriminant analysis, and so forth, the multivariate problem of testing for homogeneity of covariance matrices is certainly no less important than its univariate counterpart. The same problem moreover is of intrinsic interest in such fields as psychometrics or genetics where, for instance, the homogeneity of genetic covariance structure among species is a classical subject of investigation. Robustness and power issues, however, are even more delicate and complex in the multivariate context.

Here again, the traditional procedure is a Gaussian modified likelihood ratio test, $\phi^{(n)}_{\mathrm{MLRT}}$, the unbiasedness of which has been established (for general dimension $k$ and number $m$ of populations, under Gaussian assumptions) by Perlman [31]. Unfortunately, this test shares the poor resistance to nonnormality of its univariate counterpart, and is invalid even under elliptical densities with finite fourth-order moments: see [12, 38] or [41]. The test proposed by Nagao [27], as shown in [19], is asymptotically equivalent to $\phi^{(n)}_{\mathrm{MLRT}}$—still under finite fourth-order moments, and hence inherits of the poor robustness properties of the latter. Quite surprisingly, and except for some attempts to bootstrap $\phi^{(n)}_{\mathrm{MLRT}}$ ([11, 42, 43]), the important problem of testing homogeneity of covariance matrices under possibly non-Gaussian densities has remained an open problem for more than fifty years.

In 2001, Schott [35] proposed a Wald test, $\phi^{(n)}_{\mathrm{Schott}}$, the validity of which still requires Gaussian densities, but also two modified versions thereof—denote them as $\phi^{(n)}_{\mathrm{Schott}*}$ and $\phi^{(n)}_{\mathrm{Schott}\dagger}$—enjoying validity robustness at homokurtic and heterokurtic elliptical densities with finite fourth-order moments, respectively. More recently, a detailed account of the pseudo-Gaussian approach of the problem has been given in [19], where we derive the *locally asymptotically most stringent* (in the Le Cam sense) Gaussian test $\phi^{(n)}_{\mathcal{N}}$, and show how to turn it into pseudo-Gaussian tests $\phi^{(n)}_{\mathcal{N}*}$ and $\phi^{(n)}_{\mathcal{N}\dagger}$ that are valid under homokurtic and heterokurtic $m$-tuples of elliptical densities with finite fourth-order moments, respectively. We also show that $\phi^{(n)}_{\mathrm{Schott}*}$ and $\phi^{(n)}_{\mathcal{N}*}$ (resp. $\phi^{(n)}_{\mathrm{Schott}\dagger}$ and $\phi^{(n)}_{\mathcal{N}\dagger}$) asymptotically coincide. Schott himself recommends using the homokurtic version $\phi^{(n)}_{\mathrm{Schott}*}$ of his test and, since heterokurtic samples are not considered here, we will concentrate on $\phi^{(n)}_{\mathrm{Schott}*}$ and $\phi^{(n)}_{\mathcal{N}*}$ as the only pseudo-Gaussian tests available so far for the problem under study (one of the findings of [19] indeed is that the bootstrapped MLRT, although perfectly valid, is highly unsatisfactory in the non-Gaussian case).

This pseudo-Gaussian qualification, however, still requires finite moments of order four, and only addresses validity robustness issues. Although no multivariate equivalent of the Conover, Johnson and Johnson study [7] has been conducted so far, the results we are obtaining in Sections 5.2 and 6 below, however, indicate that $\phi^{(n)}_{\mathrm{Schott}*}$ and $\phi^{(n)}_{\mathcal{N}*}$ still suffer a severe lack of efficiency robustness, particularly so at radial densities with high kurtosis such as the Student with $4+\delta$ degrees of freedom (with $\delta \approx 0$).

It is thus very likely that the conclusions of Conover et al. [7] still apply, which strongly suggests considering a "rank-based" approach—with a concept of "ranks" adapted to the multivariate context. The purpose of this paper is to develop such an approach under elliptical assumptions, based on the signs and ranks considered in [15] for location, [16] for (auto)regression, [14] and [17] for shape.

Contrary to all existing methods, our tests do not require any moment assumptions, so that the null hypothesis they address actually is that of homogeneity of *scatter matrices*, reducing to more classical homogeneity of covariance matrices under finite second-order moments. They are asymptotically distribution-free, reach semiparametric efficiency at correctly specified densities, and are both validity- and efficiency-robust. When based on Gaussian scores, their asymptotic relative efficiency with respect to the Gaussian and pseudo-Gaussian tests is larger than one under almost all elliptical densities (see Section 6 for details).

1.2. *Testing equality of scatter (covariance) matrices.* Denote by $(\mathbf{X}_{i1}, \ldots, \mathbf{X}_{in_i})$, $i = 1, \ldots, m$, a collection of $m$ mutually independent samples of i.i.d. $k$-dimensional random vectors with *elliptically symmetric densities*. More precisely, for all $i = 1, \ldots, m$, the $n_i$ observations $\mathbf{X}_{ij}$, $j = 1, \ldots, n_i$, are assumed to have a probability density function of the form

$$\underline{f_i}(\mathbf{x}) := c_{k,f_1} |\boldsymbol{\Sigma}_i|^{-1/2} f_1(((\mathbf{x}-\boldsymbol{\theta}_i)'\boldsymbol{\Sigma}_i^{-1}(\mathbf{x}-\boldsymbol{\theta}_i))^{1/2}), \qquad \mathbf{x} \in \mathbb{R}^k, \tag{1.1}$$

for some $k$-dimensional vector $\boldsymbol{\theta}_i$ (*location*), symmetric and positive definite $(k \times k)$ matrix $\boldsymbol{\Sigma}_i$ (the *scatter* matrix), and (duly standardized: see below) function $f_1 : \mathbb{R}_0^+ \to \mathbb{R}^+$ (the *radial density*). The null hypothesis considered throughout is the hypothesis $\mathcal{H}_0 : \boldsymbol{\Sigma}_1 = \cdots = \boldsymbol{\Sigma}_m$ of *scatter homogeneity* (reducing, under finite variances, to *covariance homogeneity*).

Define (throughout, $\mathbf{\Sigma}^{1/2}$ stands for the symmetric root of $\mathbf{\Sigma}$)

$$\mathbf{U}_{ij}(\boldsymbol{\theta}_i, \mathbf{\Sigma}_i) := \frac{\mathbf{\Sigma}_i^{-1/2}(\mathbf{X}_{ij} - \boldsymbol{\theta}_i)}{\|\mathbf{\Sigma}_i^{-1/2}(\mathbf{X}_{ij} - \boldsymbol{\theta}_i)\|} \quad \text{and}$$
$$d_{ij}(\boldsymbol{\theta}_i, \mathbf{\Sigma}_i) := \|\mathbf{\Sigma}_i^{-1/2}(\mathbf{X}_{ij} - \boldsymbol{\theta}_i)\|. \tag{1.2}$$

Denoting by $n = \sum_{i=1}^m n_i$ the total sample size and writing $R_{ij}(\boldsymbol{\theta}_1, \ldots, \boldsymbol{\theta}_m, \mathbf{\Sigma}_1, \ldots, \mathbf{\Sigma}_m)$ for the rank of $d_{ij}(\boldsymbol{\theta}_i, \mathbf{\Sigma}_i)$ among $d_{11}(\boldsymbol{\theta}_1, \mathbf{\Sigma}_1), \ldots, d_{mn_m}(\boldsymbol{\theta}_m, \mathbf{\Sigma}_m)$, consider the *signed-rank scatter matrices*

$$\underset{\sim}{\mathbf{S}}_{K;i} := \frac{1}{n_i} \sum_{j=1}^{n_i} K\left(\frac{R_{ij}(\hat{\boldsymbol{\theta}}_1, \ldots, \hat{\boldsymbol{\theta}}_m, \hat{\mathbf{\Sigma}}, \ldots, \hat{\mathbf{\Sigma}})}{n+1}\right) \mathbf{U}_{ij}(\hat{\boldsymbol{\theta}}_i, \hat{\mathbf{\Sigma}}) \mathbf{U}'_{ij}(\hat{\boldsymbol{\theta}}_i, \hat{\mathbf{\Sigma}}),$$
$$i = 1, \ldots, m, \tag{1.3}$$

where $\hat{\boldsymbol{\theta}}_1, \ldots, \hat{\boldsymbol{\theta}}_m$ are consistent estimates of the various location parameters, $\hat{\mathbf{\Sigma}}$ is a consistent (under $\mathcal{H}_0$) estimate of the common null value of the $\mathbf{\Sigma}_i$'s, and $K$ is some appropriate score function. The proposed signed-rank tests reject the null hypothesis of scatter homogeneity for large values of

$$\underset{\sim}{Q}_K^{(n)} := \frac{1}{n} \sum_{1<i<i'<m} (n_i + n_{i'}) \underset{\sim}{Q}_{K;i,i'}^{(n)}, \tag{1.4}$$

where $\underset{\sim}{Q}_{K;i,i'}^{(n)} := \frac{n_i n_{i'}}{n_i + n_{i'}} \{\alpha_{k,K} \operatorname{tr}[(\underset{\sim}{\mathbf{S}}_{K;i} - \underset{\sim}{\mathbf{S}}_{K;i'})^2] + \beta_{k,K} \operatorname{tr}^2[\underset{\sim}{\mathbf{S}}_{K;i} - \underset{\sim}{\mathbf{S}}_{K;i'}]\}$; $\alpha_{k,K}$ and $\beta_{k,K}$ are constants depending on the dimension $k$ and the score function $K$. The form of (1.4) follows from Le Cam type optimality arguments, but $\underset{\sim}{Q}_K^{(n)}$ also can be obtained by replacing traditional sample covariance matrices with the signed-rank scatter matrices (1.3) in the statistic of the pseudo-Gaussian test $\phi_{\mathcal{N}*}^{(n)}$ derived in [19].

The use of signed ranks is justified by the invariance principle: $\mathcal{H}_0$ indeed is not only invariant under affine transformations, but also under the group of (continuous monotone) radial transformations; see Section 3.2 for details. Beyond affine-invariance (all tests considered in this paper are affine-invariant), our rank tests—unlike their competitors—are also (asymptotically) invariant with respect to the groups of radial transformations; their validity robustness actually follows from this latter invariance property.

As announced, our methodology combines validity and efficiency robustness. We will show that for (essentially) any radial density $f_1$, it is possible to define a score function $K := K_{f_1}$ characterizing a signed-rank test which is *locally and asymptotically optimal* (*locally and asymptotically most stringent*, in the Le Cam sense) under radial density $f_1$. In particular, when based on Gaussian scores, our rank tests achieve the same asymptotic performances as the optimal Gaussian tests $\phi_{\mathcal{N}}^{(n)}$ at the multinormal, while enjoying the validity robustness of the pseudo-Gaussian $\phi_{\mathcal{N}*}^{(n)}$ (or $\phi_{\text{Schott}*}^{(n)}$) and even more, since no moment assumption is required. Moreover, the asymptotic relative efficiencies (AREs) of these normal-score tests are almost always larger than one with respect to their parametric competitors; see Section 6. The class of tests we are proposing thus in most cases dominates the existing parametric ones, both in terms of robustness and power.

1.3. *Outline of the paper.* The paper is organized as follows. In Section 2, we collect the main assumptions needed in the sequel. Section 3 discusses semiparametric modeling issues and their relation to group invariance. Section 4.1 states the uniform local asymptotic normality result (ULAN) on which our construction of locally and asymptotically optimal tests is based. In Section 4.2, we construct rank-based versions of the *central sequences* appearing in this ULAN result. In Section 5.1, we derive and study the proposed nonparametric (signed-rank) tests [based on (1.4)] and, in Section 5.2, for the purpose of performance comparisons, the optimal pseudo-Gaussian ones. Asymptotic relative efficiencies with respect to those pseudo-Gaussian tests are derived in Section 6. Section 7 provides some simulation results confirming the theoretical ones. Finally, the Appendix collects proofs of asymptotic linearity and other technical results.

## 2. Main assumptions.

For the sake of convenience, we are collecting here the main assumptions to be used in the sequel.

2.1. *Elliptical symmetry.* As mentioned before, we throughout assume that all populations are elliptically symmetric. More precisely, defining the collections $\mathcal{F}$ of *radial densities* and $\mathcal{F}_1$ of *standardized radial densities* as

$$\mathcal{F} := \{f > 0 \text{ a.e.} : \mu_{k-1;f} < \infty\}$$

and

$$\mathcal{F}_1 := \left\{ f_1 \in \mathcal{F} : \frac{1}{\mu_{k-1;f_1}} \int_0^1 r^{k-1} f_1(r)\, dr = \frac{1}{2} \right\},$$

respectively, where $\mu_{\ell;f} := \int_0^\infty r^\ell f(r)\, dr$, we require the following.

ASSUMPTION (A). The observations $\mathbf{X}_{ij}$, $j = 1, \ldots, n_i$, $i = 1, \ldots, m$, are mutually independent, with p.d.f. $\underline{f_i}$, $i = 1, \ldots, m$, given in (1.1), for some $f_1 \in \mathcal{F}_1$.

Clearly, for the scatter matrices $\mathbf{\Sigma}_i$ in (1.1) to be well defined, identifiability restrictions are needed. This is why we impose that $f_1 \in \mathcal{F}_1$, which implies that $d_{ij}(\boldsymbol{\theta}_i, \mathbf{\Sigma}_i)$ defined in (1.2) has median one under (1.1), and identifies $\mathbf{\Sigma}_i$ without requiring any moment assumptions (see [17] for a discussion). Note, however, that under finite second-order moments, $\mathbf{\Sigma}_i$ is proportional to the covariance $\mathbf{\Sigma}_{0i}$ of $\mathbf{X}_{ij}$, with a proportionality constant that does not depend on $i$: the hypotheses of scatter and covariance homogeneity then coincide.

Special instances of elliptical densities are the $k$-variate multinormal distributions, with standardized radial density $f_1(r) = \phi_1(r) := \exp(-a_k r^2/2)$, the $k$-variate Student distributions, with radial densities (for $\nu \in \mathbb{R}_0^+$ degrees of freedom) $f_1(r) = f_{1,\nu}^t(r) := (1 + a_{k,\nu} r^2/\nu)^{-(k+\nu)/2}$, and the $k$-variate power-exponential distributions, with radial densities of the form $f_1(r) = f_{1,\eta}^e(r) := \exp(-b_{k,\eta} r^{2\eta})$, $\eta \in \mathbb{R}_0^+$; the positive constants $a_k$, $a_{k,\nu}$ and $b_{k,\eta}$ are such that $f_1 \in \mathcal{F}_1$.

The equidensity contours associated with (1.1) are hyper-ellipsoids centered at $\boldsymbol{\theta}_i$, the shape and orientation of which are determined by the scatter matrix $\mathbf{\Sigma}_i$. The *multivariate signs* $\mathbf{U}_{ij}(\boldsymbol{\theta}_i, \mathbf{\Sigma}_i)$ and *standardized radial distances* $d_{ij}(\boldsymbol{\theta}_i, \mathbf{\Sigma}_i)$ defined in (1.2) are the (within-group) elliptical coordinates associated with those ellipsoids: the observation $\mathbf{X}_{ij}$ then decomposes into $\boldsymbol{\theta}_i + d_{ij} \mathbf{\Sigma}_i^{1/2} \mathbf{U}_{ij}$, where the $\mathbf{U}_{ij}$'s, $j = 1, \ldots, n_i$, $i = 1, \ldots, m$ are i.i.d. uniform over the unit sphere in $\mathbb{R}^k$, and the $d_{ij}$'s are i.i.d., independent of the $\mathbf{U}_{ij}$'s, with common density $\tilde{f}_{1k}(r) := (\mu_{k-1;f_1})^{-1} r^{k-1} f_1(r)$ over $\mathbb{R}^+$ (justifying the terminology *standardized radial density* for $f_1$) and distribution function $\tilde{F}_{1k}$. In the sequel, the notation $\tilde{g}_{1k}$ and $\tilde{G}_{1k}$ will be used for the corresponding functions computed from a standardized radial density $g_1 (\in \mathcal{F}_1)$.

The derivation of locally and asymptotically optimal tests at radial density $f_1$ will be based on the *uniform local and asymptotic normality* (ULAN) of the model *at given* $f_1$. This ULAN property—the statement of which requires some further preparation and is delayed to Section 4.1—only holds under some further mild regularity conditions on $f_1$. More precisely, ULAN (see Proposition 4.1 below) requires $f_1$ to belong to the collection $\mathcal{F}_a$ of absolutely continuous densities in $\mathcal{F}_1$ such that, letting $\varphi_{f_1} := -\dot{f}_1/f_1$ (with $\dot{f}_1$ the a.e.-derivative of $f_1$), the integrals

$$\mathcal{I}_k(f_1) := \int_0^1 \varphi_{f_1}^2(\tilde{F}_{1k}^{-1}(u))\, du \quad \text{and} \quad \mathcal{J}_k(f_1) := \int_0^1 \varphi_{f_1}^2(\tilde{F}_{1k}^{-1}(u))(\tilde{F}_{1k}^{-1}(u))^2\, du$$

are finite. The quantities $\mathcal{I}_k(f_1)$ and $\mathcal{J}_k(f_1)$ play the roles of *radial Fisher information for location* and *radial Fisher information for shape/scale*, respectively (see [17]).

2.2. *Asymptotic behavior of sample sizes.* Although, for the sake of notational simplicity, we do not mention it explicitly, we actually consider sequences of statistical experiments, with triangular arrays of observations of the form $(\mathbf{X}_{11}^{(n)}, \ldots, \mathbf{X}_{1n_1^{(n)}}^{(n)}, \mathbf{X}_{21}^{(n)}, \ldots, \mathbf{X}_{2n_2^{(n)}}^{(n)}, \ldots, \mathbf{X}_{m1}^{(n)}, \ldots, \mathbf{X}_{mn_m^{(n)}}^{(n)})$ indexed by the total sample size $n$, where the sequences $n_i^{(n)}$ satisfy the following assumption.

ASSUMPTION (B). For all $i = 1, \ldots, m$, $n_i = n_i^{(n)} \to \infty$ as $n \to \infty$.

Note that this assumption is weaker than the corresponding classical assumption in (univariate or multivariate) multisample problems, which requires that $n_i/n$ be bounded away from 0 and 1 for all $i$ as $n \to \infty$. Letting $\lambda_i^{(n)} := n_i^{(n)}/n$, it is easy to check that Assumption (B) is actually equivalent to the Noether conditions

$$\max\left(\frac{1-\lambda_i^{(n)}}{\lambda_i^{(n)}}, \frac{\lambda_i^{(n)}}{1-\lambda_i^{(n)}}\right) = o(n) \qquad \text{as } n \to \infty, \text{ for all } i,$$

that are needed for the representation result in Lemma 4.1(i) below. However, in the derivation of asymptotic distributions under local alternatives, the following reinforcement of Assumption (B) is assumed to hold (mainly, for notational comfort):

ASSUMPTION (B′). For all $i = 1, \ldots, m$, $\lambda_i^{(n)} \to \lambda_i \in (0, 1)$, as $n \to \infty$.

2.3. *Score functions.* The score functions $K$ appearing in the rank-based statistics (1.3)–(1.4) will be assumed to satisfy the following regularity conditions.

ASSUMPTION (C). The score function $K : (0, 1) \to \mathbb{R}$ (C1) is a continuous, nonconstant, and square-integrable mapping which (C2) can be expressed as the difference of two monotone increasing functions, and (C3) satisfies $\int_0^1 K(u)\,du = k$.

Assumption (C3) is a normalization constraint that is automatically satisfied by the score functions $K(u) = K_{f_1}(u) := \varphi_{f_1}(\tilde{F}_{1k}^{-1}(u))\tilde{F}_{1k}^{-1}(u)$ leading to local and asymptotic optimality at radial density $f_1$ (at which ULAN holds); see Section 4.1. For score functions $K, K_1, K_2$ satisfying Assumption (C), let $\mathcal{J}_k(K_1, K_2) := \mathrm{E}[K_1(U)K_2(U)]$ and $\mathcal{L}_k(K_1, K_2) := \mathrm{Cov}[K_1(U), K_2(U)] = \mathcal{J}_k(K_1, K_2) - k^2$ [throughout, $U$ stands for a random variable uniformly distributed over $(0, 1)$]. For simplicity, we write $\mathcal{J}_k(K)$ for $\mathcal{J}_k(K, K)$, $\mathcal{L}_k(K)$ for

$\mathcal{L}_k(K,K)$, $\mathcal{J}_k(K,f_1)$ for $\mathrm{E}[K(U)K_{f_1}(U)]$, $\mathcal{L}_k(f_1,g_1)$ for $\mathrm{E}[K_{f_1}(U)K_{g_1}(U)]-k^2$, etc.

The power score functions $K_a(u):=k(a+1)u^a$ $(a>0)$ provide some traditional score functions satisfying Assumption (C), with $\mathcal{J}_k(K_a)=k^2(a+1)^2/(2a+1)$ and $\mathcal{L}_k(K_a)=k^2a^2/(2a+1)$: Wilcoxon and Spearman scores are obtained for $a=1$ and $a=2$, respectively. As for score functions of the form $K_{f_1}$, an important particular case is that of van der Waerden or normal scores, obtained for $f_1=\phi_1$. Then, denoting by $\Psi_k$ the chi-square distribution function with $k$ degrees of freedom, $K_{\phi_1}(u)=\Psi_k^{-1}(u)$, $\mathcal{J}_k(\phi_1)=k(k+2)$, and $\mathcal{L}_k(\phi_1)=2k$. Similarly, writing $G_{k,\nu}$ for the Fisher–Snedecor distribution function with $k$ and $\nu$ degrees of freedom, Student densities $f_1=f_{1,\nu}^t$ yield $K_{f_{1,\nu}^t}(u)=k(k+\nu)G_{k,\nu}^{-1}(u)/(\nu+kG_{k,\nu}^{-1}(u))$, $\mathcal{J}_k(f_{1,\nu}^t)=k(k+2)(k+\nu)/(k+\nu+2)$ and $\mathcal{L}_k(f_{1,\nu}^t)=2k\nu/(k+\nu+2)$.

## 3. Semiparametric modeling of the family of elliptical densities.

3.1. *Scatter, scale, and shape.* Consider an observed $n$-tuple $\mathbf{X}_1,\dots,\mathbf{X}_n$ of i.i.d. $k$-dimensional elliptical random vectors, with location $\boldsymbol{\theta}$, scatter $\boldsymbol{\Sigma}$, and radial density $f_1\in\mathcal{F}_1$ but otherwise unspecified. The family $\mathcal{P}^{(n)}$ of distributions for this observation is indexed by $(\boldsymbol{\theta},\boldsymbol{\Sigma},f_1)$. As soon as a semiparametric point of view is adopted, or when rank-based methods are considered, the scatter matrix $\boldsymbol{\Sigma}$ naturally decomposes into $\boldsymbol{\Sigma}=\sigma^2\mathbf{V}$, where $\sigma$ is a *scale parameter* (equivariant under multiplication by a positive constant) and $\mathbf{V}$ a *shape matrix* (invariant under multiplication by a positive constant). A semiparametric model with specified $\sigma$ and unspecified standardized radial density $f_1$ indeed would be highly artificial, and we, therefore, only consider the case under which $\sigma$ and $f_1$ are both unspecified. This semiparametric setting is also the one that enjoys the group invariance structure in which the ranks and the signs to be used in our method spontaneously arise from invariance arguments; see Section 3.2 below.

The concepts of scale and shape however require a more careful definition. Denoting by $\mathcal{S}_k$ the collection of all $k\times k$ symmetric positive definite real matrices, consider a function $S:\mathcal{S}_k\to\mathbb{R}_0^+$ satisfying $S(\lambda\boldsymbol{\Sigma})=\lambda S(\boldsymbol{\Sigma})$ for all $\lambda\in\mathbb{R}_0^+$ and $\boldsymbol{\Sigma}\in\mathcal{S}_k$, and define scale and shape as $\sigma_S:=(S(\boldsymbol{\Sigma}))^{1/2}$ and $\mathbf{V}_S:=\boldsymbol{\Sigma}/S(\boldsymbol{\Sigma})$, respectively. Clearly, $\mathbf{V}_S$ is the only matrix in $\mathcal{S}_k$ which is proportional to $\boldsymbol{\Sigma}$ and satisfies $S(\mathbf{V}_S)=1$: denote by $\mathcal{V}_k^S:=\{\mathbf{V}\in\mathcal{S}_k : S(\mathbf{V})=1\}$ the set of all possible shape matrices associated with $S$. Classical choices of $S$ are (i) $S(\boldsymbol{\Sigma})=(\boldsymbol{\Sigma})_{11}$ ([14, 17, 23] and [33]); (ii) $S(\boldsymbol{\Sigma})=k^{-1}\operatorname{tr}(\boldsymbol{\Sigma})$ ([8, 28] and [39]); (iii) $S(\boldsymbol{\Sigma})=|\boldsymbol{\Sigma}|^{1/k}$ ([9, 34, 36] and [37]).

In practice, all choices of $S$ are essentially equivalent. Bickel (Example 4 of [4], with a trace-based normalization of $\boldsymbol{\Sigma}^{-1}$) actually shows that irrespective of $S$, the asymptotic information matrix for $\mathbf{V}_S$ in the presence of unspecified $\boldsymbol{\theta}$ and $\sigma_S$ is the same, at any $\boldsymbol{\theta} \in \mathbb{R}^k$, $\sigma_S \in \mathbb{R}_0^+$, $\mathbf{V}_S \in \mathcal{V}_k^S$ and $f_1$, whether $f_1$ is specified (parametric model) or not (semiparametric model): once $\boldsymbol{\theta}$ and $\sigma_S$ are unspecified, an unspecified $f_1$ does not induce any additional loss for inference about $\mathbf{V}_S$. In [30], Paindaveine establishes the stronger result that the information matrix for $\mathbf{V}_S$ in the presence of unspecified $\boldsymbol{\theta}$, $\sigma_S$ and $f_1$ is strictly less, at any $\boldsymbol{\theta} \in \mathbb{R}^k$, $\sigma_S \in \mathbb{R}_0^+$, $\mathbf{V}_S \in \mathcal{V}_k^S$ and $f_1$, than in the corresponding parametric model with specified $\boldsymbol{\theta}$, $\sigma_S$ and $f_1$—except for $S : \boldsymbol{\Sigma} \mapsto |\boldsymbol{\Sigma}|^{1/k}$, where those two information matrices coincide: under this determinant-based normalization, thus, the presence of nuisances ($\boldsymbol{\theta}$, $\sigma_S$ and $f_1$) (resp., $\boldsymbol{\theta}$, $\mathbf{V}_S$, and $f_1$) asymptotically has no effect on inference about shape (resp., inference about scale). In both cases, it can be said (adopting a point estimation terminology) that shape can be estimated *adaptively*. This *Paindaveine adaptivity*, however, where $\boldsymbol{\theta}$, $\sigma_S$ and $f_1$ lie in the nuisance space of the semiparametric model, is much stronger than *Bickel adaptivity* where only $f_1$ does. This finding strongly pleads in favor of the determinant-based definition of shape which, with its block-diagonal information matrix for $\boldsymbol{\theta}, \sigma_S$ and $\mathbf{V}_S$, is also more convenient from the point of view of statistical inference, and as we will see in Section 5.1, allows for an ANOVA-type decomposition of the test statistics into two mutually independent parts providing tests against subalternatives of scale and shape heterogeneity, respectively. Therefore, throughout we adopt $S(\boldsymbol{\Sigma}) = |\boldsymbol{\Sigma}|^{1/k}$, and henceforth simply write $\mathbf{V} \in \mathcal{V}_k$ and $\sigma$ for the resulting shape and scale.

The parameter in our problem then is the $L$-dimensional vector

$$\boldsymbol{\vartheta} := (\boldsymbol{\vartheta}_I', \boldsymbol{\vartheta}_{II}', \boldsymbol{\vartheta}_{III}')' := (\boldsymbol{\theta}_1', \ldots, \boldsymbol{\theta}_m', \sigma_1^2, \ldots, \sigma_m^2, (\overset{\circ}{\mathrm{vech}} \mathbf{V}_1)', \ldots, (\overset{\circ}{\mathrm{vech}}\, \mathbf{V}_m)')',$$

where $L = mk(k+3)/2$ and $\overset{\circ}{\mathrm{vech}}(\mathbf{V})$ is characterized by $\mathrm{vech}(\mathbf{V}) =: ((\mathbf{V})_{11}, (\overset{\circ}{\mathrm{vech}}\, \mathbf{V})')'$. Indeed, $\boldsymbol{\Sigma}_i$ is entirely determined by $\sigma_i^2$ and $\overset{\circ}{\mathrm{vech}}(\mathbf{V}_i)$. Write $\boldsymbol{\Theta}$ for the set $\mathbb{R}^{mk} \times (\mathbb{R}_0^+)^m \times (\overset{\circ}{\mathrm{vech}}(\mathcal{V}_k))^m$ of admissible $\boldsymbol{\vartheta}$ values, and $\mathrm{P}_{\boldsymbol{\vartheta}; f_1}^{(n)}$ or $\mathrm{P}_{\boldsymbol{\vartheta}_I, \boldsymbol{\vartheta}_{II}, \boldsymbol{\vartheta}_{III}; f_1}^{(n)}$ for the joint distribution of the $n$ observations under parameter value $\boldsymbol{\vartheta}$ and standardized radial density $f_1$ (always implicitly assumed to belong to $\mathcal{F}_1$, when notation $f_1$ is used).

Finally, note that (i) $\mathbf{U}_{ij}(\boldsymbol{\theta}_i, \boldsymbol{\Sigma}_i) = \mathbf{U}_{ij}(\boldsymbol{\theta}_i, \mathbf{V}_i)$ and (ii) $d_{ij}(\boldsymbol{\theta}_i, \boldsymbol{\Sigma}_i) = \sigma_i^{-1} d_{ij}(\boldsymbol{\theta}_i, \mathbf{V}_i)$. It follows from (ii) that under $\mathcal{H}_0$, the ranks of the radial distances computed from the common value $\mathbf{V}$ of the shape matrices coincide with those of the standardized radial distances computed from the common value $\boldsymbol{\Sigma}$ of the scatter matrices.

3.2. *Invariance issues.* Denoting by $\mathcal{M}(\boldsymbol{\Upsilon})$ the vector space spanned by the columns of the $L \times r$ full-rank matrix $\boldsymbol{\Upsilon}$ $(r < L)$, the null hypothesis of scatter homogeneity $\mathcal{H}_0 : \sigma_1^2 \mathbf{V}_1 = \cdots = \sigma_m^2 \mathbf{V}_m$ can be written as $\mathcal{H}_0 : \boldsymbol{\vartheta} \in \mathcal{M}(\boldsymbol{\Upsilon})$, with

$$\boldsymbol{\Upsilon} := \begin{pmatrix} \boldsymbol{\Upsilon}_I & \mathbf{0} & \mathbf{0} \\ \mathbf{0} & \boldsymbol{\Upsilon}_{II} & \mathbf{0} \\ \mathbf{0} & \mathbf{0} & \boldsymbol{\Upsilon}_{III} \end{pmatrix} := \begin{pmatrix} \mathbf{I}_{mk} & \mathbf{0} & \mathbf{0} \\ \mathbf{0} & \mathbf{1}_m & \mathbf{0} \\ \mathbf{0} & \mathbf{0} & \mathbf{1}_m \otimes \mathbf{I}_{k_0} \end{pmatrix}, \tag{3.1}$$
$$k_0 := \frac{k(k+1)}{2} - 1,$$

where $\mathbf{1}_m := (1, \ldots, 1)' \in \mathbb{R}^m$ and $\mathbf{I}_\ell$ denotes the $\ell$-dimensional identity matrix.

Two distinct invariance structures play a role here. The first one is related with the group of affine transformations of the observations, which generates the parametric families $\mathcal{P}^{(n)}_{\boldsymbol{\Upsilon}, f_1} := \bigcup_{\boldsymbol{\vartheta} \in \mathcal{M}(\boldsymbol{\Upsilon})} \{\mathrm{P}^{(n)}_{\boldsymbol{\vartheta}; f_1}\}$. More precisely, this group is the group $\mathcal{G}^{m,k}, \circ$ of affine transformations of the form $\mathbf{X}_{ij} \mapsto \mathbf{A}\mathbf{X}_{ij} + \mathbf{b}_i$, where $\mathbf{A}$ is a full-rank $(k \times k)$ matrix and $\mathbf{B} := (\mathbf{b}_1, \ldots, \mathbf{b}_m)$ a $(k \times m)$ matrix. Associated with that group is the group $\tilde{\mathcal{G}}^{m,k}, \circ$ of transformations $\boldsymbol{\vartheta} \mapsto \mathbf{g}^{m,k}_{\mathbf{A},\mathbf{B}}(\boldsymbol{\vartheta})$ of the parameter space, where

$$\begin{aligned} \mathbf{g}^{m,k}_{\mathbf{A},\mathbf{B}}(\boldsymbol{\vartheta}) &:= ((\mathbf{A}\boldsymbol{\theta}_1 + \mathbf{b}_1)', \ldots, (\mathbf{A}\boldsymbol{\theta}_m + \mathbf{b}_m)', |\mathbf{A}|^{2/k}\sigma_1^2, \ldots, |\mathbf{A}|^{2/k}\sigma_m^2, \\ &\qquad (\overset{\circ}{\mathrm{vech}}(\mathbf{A}\mathbf{V}_1\mathbf{A}'))'/|\mathbf{A}|^{2/k}, \ldots, (\overset{\circ}{\mathrm{vech}}(\mathbf{A}\mathbf{V}_m\mathbf{A}'))'/|\mathbf{A}|^{2/k})'. \end{aligned}$$

Clearly, $\mathcal{H}_0$ is invariant under $\mathcal{G}^{m,k}, \circ$—meaning that $\mathbf{g}^{m,k}_{\mathbf{A},\mathbf{B}}(\mathcal{M}(\boldsymbol{\Upsilon})) = \mathcal{M}(\boldsymbol{\Upsilon})$ for all $\mathbf{g}^{m,k}_{\mathbf{A},\mathbf{B}}$. Therefore, it is reasonable to restrict to affine-invariant tests of $\mathcal{H}_0$. Beyond their distribution-freeness with respect to the $\boldsymbol{\theta}_i$'s and the common null values $\sigma$ and $\mathbf{V}$ of the scale and shape parameters, affine-invariant test statistics—that is, statistics $Q$ such that $Q(\mathbf{A}\mathbf{X}_{11} + \mathbf{b}_1, \ldots, \mathbf{A}\mathbf{X}_{mn_m} + \mathbf{b}_m) = Q(\mathbf{X}_{11}, \ldots, \mathbf{X}_{mn_m})$ for all $\mathbf{A}, \mathbf{b}_1, \ldots, \mathbf{b}_m$—yield tests that are *coordinate-free.*

A second invariance structure is induced by the groups $\mathcal{G}, \circ := \mathcal{G}^{\boldsymbol{\vartheta}_I, \mathbf{V}}, \circ$ of *continuous monotone radial transformations*, of the form

$$\begin{aligned} \mathbf{X} &\mapsto \mathcal{G}_g(\mathbf{X}_{11}, \ldots, \mathbf{X}_{mn_m}) \\ &= \mathcal{G}_g(\boldsymbol{\theta}_1 + d_{11}(\boldsymbol{\theta}_1, \mathbf{V})\mathbf{V}^{1/2}\mathbf{U}_{11}(\boldsymbol{\theta}_1, \mathbf{V}), \ldots, \\ &\qquad \boldsymbol{\theta}_m + d_{mn_m}(\boldsymbol{\theta}_m, \mathbf{V})\mathbf{V}^{1/2}\mathbf{U}_{mn_m}(\boldsymbol{\theta}_m, \mathbf{V})) \\ &:= (\boldsymbol{\theta}_1 + g(d_{11}(\boldsymbol{\theta}_1, \mathbf{V}))\mathbf{V}^{1/2}\mathbf{U}_{11}(\boldsymbol{\theta}_1, \mathbf{V}), \ldots, \\ &\qquad \boldsymbol{\theta}_m + g(d_{mn_m}(\boldsymbol{\theta}_m, \mathbf{V}))\mathbf{V}^{1/2}\mathbf{U}_{mn_m}(\boldsymbol{\theta}_m, \mathbf{V})), \end{aligned}$$

where $g : \mathbb{R}^+ \to \mathbb{R}^+$ is continuous, monotone increasing, and such that $g(0) = 0$ and $\lim_{r \to \infty} g(r) = \infty$. For each $\boldsymbol{\vartheta}_I, \mathbf{V}$, this group $\mathcal{G}^{\boldsymbol{\vartheta}_I, \mathbf{V}}, \circ$ is a generating group for $\mathcal{P}^{(n)}_{\boldsymbol{\vartheta}_I,\mathbf{V}} := \bigcup_\sigma \bigcup_{f_1} \{\mathrm{P}^{(n)}_{\boldsymbol{\vartheta}_I,\sigma^2\mathbf{1}_m,\mathbf{1}_m\otimes(\overset{\circ}{\mathrm{vech}}\mathbf{V});f_1}\}$ (a nonparametric family). In such families, the invariance principle suggests basing inference on statistics that are measurable with respect to the corresponding maximal invariant, namely the vectors $(\mathbf{U}_{11},\dots,\mathbf{U}_{mn_m})$ of signs and the vector $(R_{11},\dots,R_{mn_m})$ of ranks, where $\mathbf{U}_{ij}=\mathbf{U}_{ij}(\boldsymbol{\theta}_i,\mathbf{V})$, and $R_{ij}=R_{ij}(\boldsymbol{\theta}_1,\dots,\boldsymbol{\theta}_m,\mathbf{V},\dots,\mathbf{V})$. Such invariant statistics of course are distribution-free under $\mathcal{P}^{(n)}_{\boldsymbol{\vartheta}_I,\mathbf{V}}$.

## 4. Uniform local asymptotic normality, signs and ranks.

4.1. *Uniform local asymptotic normality* (*ULAN*). As mentioned in Section 1, we plan to develop tests that are optimal at correctly specified densities, in the sense of Le Cam's asymptotic theory of statistical experiments. In this section, we state the ULAN result (with respect to location, scale and shape parameters for fixed radial density $f_1$) on which optimality will be based. Writing

$$\begin{aligned}\boldsymbol{\vartheta}^{(n)} &= (\boldsymbol{\vartheta}_I^{(n)\prime},\boldsymbol{\vartheta}_{II}^{(n)\prime},\boldsymbol{\vartheta}_{III}^{(n)\prime})' \\ &= (\boldsymbol{\theta}_1^{(n)\prime},\dots,\boldsymbol{\theta}_m^{(n)\prime},\sigma_1^{2(n)},\dots,\sigma_m^{2(n)},(\overset{\circ}{\mathrm{vech}}\mathbf{V}_1^{(n)})',\dots,(\overset{\circ}{\mathrm{vech}}\mathbf{V}_m^{(n)})')'\end{aligned}$$

for an arbitrary sequence of $L$-dimensional parameter values in $\boldsymbol{\Theta}$, consider sequences of "local alternatives" $\boldsymbol{\vartheta}^{(n)}+n^{-1/2}\boldsymbol{\nu}^{(n)}\boldsymbol{\tau}^{(n)}$, where

$$\begin{aligned}\boldsymbol{\tau}^{(n)} &= (\boldsymbol{\tau}_I^{(n)\prime},\boldsymbol{\tau}_{II}^{(n)\prime},\boldsymbol{\tau}_{III}^{(n)\prime})' \\ &= (\mathbf{t}_1^{(n)\prime},\dots,\mathbf{t}_m^{(n)\prime},s_1^{2(n)},\dots,s_m^{2(n)},(\overset{\circ}{\mathrm{vech}}\,\mathbf{v}_1^{(n)})',\dots,(\overset{\circ}{\mathrm{vech}}\,\mathbf{v}_m^{(n)})')'\end{aligned}$$

is such that $\sup_n \boldsymbol{\tau}^{(n)\prime}\boldsymbol{\tau}^{(n)}<\infty$ and where, denoting by $\boldsymbol{\Lambda}^{(n)}=(\Lambda_{ii'}^{(n)})$ the $(m\times m)$ diagonal matrix with $\Lambda_{ii}^{(n)}:=(\lambda_i^{(n)})^{-1/2}$ (see Section 2.2),

$$\boldsymbol{\nu}^{(n)} := \begin{pmatrix}\boldsymbol{\nu}_I^{(n)} & \mathbf{0} & \mathbf{0}\\ \mathbf{0} & \boldsymbol{\nu}_{II}^{(n)} & \mathbf{0}\\ \mathbf{0} & \mathbf{0} & \boldsymbol{\nu}_{III}^{(n)}\end{pmatrix} := \begin{pmatrix}\boldsymbol{\Lambda}^{(n)}\otimes\mathbf{I}_k & \mathbf{0} & \mathbf{0}\\ \mathbf{0} & \boldsymbol{\Lambda}^{(n)} & \mathbf{0}\\ \mathbf{0} & \mathbf{0} & \boldsymbol{\Lambda}^{(n)}\otimes\mathbf{I}_{k_0}\end{pmatrix} \tag{4.1}$$

[under Assumption (B′), we also write $\boldsymbol{\nu}$ for $\lim_{n\to\infty}\boldsymbol{\nu}^{(n)}$]. Clearly, these local alternatives do not involve $(\mathbf{v}_i^{(n)})_{11}$, $i=1,\dots,m$. It is natural, though, to see that the perturbed shapes $\mathbf{V}_i^{(n)}+n_i^{-1/2}\mathbf{v}_i^{(n)}$ remain [up to $o(n_i^{-1/2})$'s] within the family $\mathcal{V}_k$ of shape matrices: this leads to defining $(\mathbf{v}_i^{(n)})_{11}$ in such a way that $\mathrm{tr}[(\mathbf{V}_i^{(n)})^{-1}\mathbf{v}_i^{(n)}]=0$, $i=1,\dots,m$, which entails $|\mathbf{V}_i^{(n)}+n_i^{-1/2}\mathbf{v}_i^{(n)}|^{1/k}=1+o(n_i^{-1/2})$ (see [18], Section 4).

The following notation will be used throughout. Let $\mathrm{diag}(\mathbf{B}_1,\mathbf{B}_2,\dots,\mathbf{B}_m)$ stand for the block-diagonal matrix with diagonal blocks $\mathbf{B}_1$, $\mathbf{B}_2,\dots,\mathbf{B}_m$.

Write $\mathbf{V}^{\otimes 2}$ for the Kronecker product $\mathbf{V} \otimes \mathbf{V}$. Denoting by $\mathbf{e}_\ell$ the $\ell$th vector of the canonical basis of $\mathbb{R}^k$, let $\mathbf{K}_k := \sum_{i,j=1}^k (\mathbf{e}_i \mathbf{e}_j') \otimes (\mathbf{e}_j \mathbf{e}_i')$ be the $(k^2 \times k^2)$ *commutation matrix*, and put $\mathbf{J}_k := (\operatorname{vec} \mathbf{I}_k)(\operatorname{vec} \mathbf{I}_k)'$. Finally, let $\mathbf{M}_k(\mathbf{V})$ be the $(k_0 \times k^2)$ matrix such that $(\mathbf{M}_k(\mathbf{V}))'(\overset{\circ}{\operatorname{vech}} \mathbf{v}) = (\operatorname{vec} \mathbf{v})$ for any symmetric $k \times k$ matrix $\mathbf{v}$ such that $\operatorname{tr}(\mathbf{V}^{-1}\mathbf{v}) = 0$. As shown in Lemma 4.2(v) of [30], $\mathbf{M}_k(\mathbf{V})(\operatorname{vec} \mathbf{V}^{-1}) = \mathbf{0}$ for all $\mathbf{V} \in \mathcal{V}_k$.

We then have the following ULAN result; the proof follows along the same lines as in Theorem 2.1 of [30], and hence is omitted.

PROPOSITION 4.1. *Assume that* (A) *and* (B) *hold, and that* $f_1 \in \mathcal{F}_a$. *Then the family* $\mathcal{P}_{f_1}^{(n)} := \{\mathrm{P}_{\boldsymbol{\vartheta}; f_1}^{(n)} | \boldsymbol{\vartheta} \in \boldsymbol{\Theta}\}$ *is ULAN, with central sequence*

$$\boldsymbol{\Delta}_{\boldsymbol{\vartheta}; f_1} = \boldsymbol{\Delta}_{\boldsymbol{\vartheta}; f_1}^{(n)} := \begin{pmatrix} \boldsymbol{\Delta}_{\boldsymbol{\vartheta}; f_1}^{I} \\ \boldsymbol{\Delta}_{\boldsymbol{\vartheta}; f_1}^{II} \\ \boldsymbol{\Delta}_{\boldsymbol{\vartheta}; f_1}^{III} \end{pmatrix},$$

$$\boldsymbol{\Delta}_{\boldsymbol{\vartheta}; f_1}^{I} = \begin{pmatrix} \boldsymbol{\Delta}_{\boldsymbol{\vartheta}; f_1}^{I,1} \\ \vdots \\ \boldsymbol{\Delta}_{\boldsymbol{\vartheta}; f_1}^{I,m} \end{pmatrix},$$

$$\boldsymbol{\Delta}_{\boldsymbol{\vartheta}; f_1}^{II} = \begin{pmatrix} \Delta_{\boldsymbol{\vartheta}; f_1}^{II,1} \\ \vdots \\ \Delta_{\boldsymbol{\vartheta}; f_1}^{II,m} \end{pmatrix}, \qquad \boldsymbol{\Delta}_{\boldsymbol{\vartheta}; f_1}^{III} = \begin{pmatrix} \boldsymbol{\Delta}_{\boldsymbol{\vartheta}; f_1}^{III,1} \\ \vdots \\ \boldsymbol{\Delta}_{\boldsymbol{\vartheta}; f_1}^{III,m} \end{pmatrix},$$

*where [with* $d_{ij} = d_{ij}(\boldsymbol{\theta}_i, \mathbf{V}_i)$ *and* $\mathbf{U}_{ij} = \mathbf{U}_{ij}(\boldsymbol{\theta}_i, \mathbf{V}_i)$*]*

$$\boldsymbol{\Delta}_{\boldsymbol{\vartheta}; f_1}^{I,i} := \frac{n_i^{-1/2}}{\sigma_i} \sum_{j=1}^{n_i} \varphi_{f_1}\left(\frac{d_{ij}}{\sigma_i}\right) \mathbf{V}_i^{-1/2} \mathbf{U}_{ij},$$

$$\Delta_{\boldsymbol{\vartheta}; f_1}^{II,i} := \frac{n_i^{-1/2}}{2\sigma_i^2} \sum_{j=1}^{n_i} \left( \varphi_{f_1}\left(\frac{d_{ij}}{\sigma_i}\right) \frac{d_{ij}}{\sigma_i} - k \right),$$

$$\boldsymbol{\Delta}_{\boldsymbol{\vartheta}; f_1}^{III,i} := \frac{n_i^{-1/2}}{2} \mathbf{M}_k(\mathbf{V}_i)(\mathbf{V}_i^{\otimes 2})^{-1/2} \sum_{j=1}^{n_i} \varphi_{f_1}\left(\frac{d_{ij}}{\sigma_i}\right) \frac{d_{ij}}{\sigma_i} \operatorname{vec}(\mathbf{U}_{ij} \mathbf{U}_{ij}'),$$

$i = 1, \ldots, m$, *and full-rank block-diagonal information matrix*

$$\boldsymbol{\Gamma}_{\boldsymbol{\vartheta}; f_1} := \operatorname{diag}(\boldsymbol{\Gamma}_{\boldsymbol{\vartheta}; f_1}^{I}, \boldsymbol{\Gamma}_{\boldsymbol{\vartheta}; f_1}^{II}, \boldsymbol{\Gamma}_{\boldsymbol{\vartheta}; f_1}^{III}), \tag{4.2}$$

*where, defining* $\underline{\boldsymbol{\sigma}} := \operatorname{diag}(\sigma_1, \ldots, \sigma_m)$, $\underline{\mathbf{V}} := \operatorname{diag}(\mathbf{V}_1, \ldots, \mathbf{V}_m)$, $\mathbf{M}_k(\underline{\mathbf{V}}) := \operatorname{diag}(\mathbf{M}_k(\mathbf{V}_1), \ldots, \mathbf{M}_k(\mathbf{V}_m))$ *and* $\underline{\mathbf{V}^{\otimes 2}} := \operatorname{diag}(\mathbf{V}_1^{\otimes 2}, \ldots, \mathbf{V}_m^{\otimes 2})$, *we let*

$$\boldsymbol{\Gamma}_{\boldsymbol{\vartheta}; f_1}^{I} := \frac{1}{k} \mathcal{I}_k(f_1)(\underline{\boldsymbol{\sigma}}^{-2} \otimes \mathbf{I}_k)\underline{\mathbf{V}}^{-1}, \qquad \boldsymbol{\Gamma}_{\boldsymbol{\vartheta}; f_1}^{II} := \tfrac{1}{4} \mathcal{L}_k(f_1)\underline{\boldsymbol{\sigma}}^{-4}$$

*and*

$$\mathbf{\Gamma}^{III}_{\boldsymbol{\vartheta};f_1} := \frac{\mathcal{J}_k(f_1)}{4k(k+2)}\mathbf{M}_k(\underline{\mathbf{V}})[\mathbf{I}_m \otimes (\mathbf{I}_{k^2} + \mathbf{K}_k)](\underline{\mathbf{V}^{\otimes 2}})^{-1}(\mathbf{M}_k(\underline{\mathbf{V}}))'.$$

*More precisely, for any* $\boldsymbol{\vartheta}^{(n)} = \boldsymbol{\vartheta} + O(n^{-1/2})$ *and any bounded sequence* $\boldsymbol{\tau}^{(n)}$*, we have*

$$\begin{aligned}\Lambda^{(n)}_{\boldsymbol{\vartheta}^{(n)}+n^{-1/2}\boldsymbol{\nu}^{(n)}\boldsymbol{\tau}^{(n)}/\boldsymbol{\vartheta}^{(n)};f_1} &:= \log(d\mathrm{P}^{(n)}_{\boldsymbol{\vartheta}^{(n)}+n^{-1/2}\boldsymbol{\nu}^{(n)}\boldsymbol{\tau}^{(n)};f_1}/d\mathrm{P}^{(n)}_{\boldsymbol{\vartheta}^{(n)};f_1}) \\ &= (\boldsymbol{\tau}^{(n)})'\mathbf{\Delta}^{(n)}_{\boldsymbol{\vartheta}^{(n)};f_1} - \tfrac{1}{2}(\boldsymbol{\tau}^{(n)})'\mathbf{\Gamma}_{\boldsymbol{\vartheta};f_1}\boldsymbol{\tau}^{(n)} + o_{\mathrm{P}}(1)\end{aligned}$$

*and* $\mathbf{\Delta}_{\boldsymbol{\vartheta}^{(n)};f_1} \xrightarrow{\mathcal{L}} \mathcal{N}(\mathbf{0}, \mathbf{\Gamma}_{\boldsymbol{\vartheta};f_1})$ *under* $\mathrm{P}^{(n)}_{\boldsymbol{\vartheta}^{(n)};f_1}$*, as* $n \to \infty$*.*

The classical theory of hypothesis testing in Gaussian shifts (see Section 11.9 of [26]) provides the general form for locally asymptotically optimal (namely, *most stringent*) tests in ULAN models. Such tests, for a null hypothesis of the form $\boldsymbol{\vartheta} \in \mathcal{M}(\mathbf{\Upsilon})$, should be based on the asymptotically chi-square null distribution of

$$Q_{\mathbf{\Upsilon}} := \mathbf{\Delta}'_{\boldsymbol{\vartheta};f_1}\mathbf{\Gamma}^{-1/2}_{\boldsymbol{\vartheta};f_1}[\mathbf{I} - \mathrm{proj}(\mathbf{\Gamma}^{1/2}_{\boldsymbol{\vartheta};f_1}(\boldsymbol{\nu}^{(n)})^{-1}\mathbf{\Upsilon})]\mathbf{\Gamma}^{-1/2}_{\boldsymbol{\vartheta};f_1}\mathbf{\Delta}_{\boldsymbol{\vartheta};f_1}$$

[with $\boldsymbol{\vartheta}$ replaced by an appropriate estimator $\hat{\boldsymbol{\vartheta}}$; see Assumption (D) below], where $\mathrm{proj}(\mathbf{A}) = \mathbf{A}[\mathbf{A}'\mathbf{A}]^{-1}\mathbf{A}'$, for any $(L \times r)$ matrix $\mathbf{A}$ with full rank $r$, is the matrix projecting $\mathbb{R}^L$ onto $\mathcal{M}(\mathbf{A})$. Whenever $\mathbf{\Gamma}_{\boldsymbol{\vartheta};f_1}$, $\boldsymbol{\nu}^{(n)}$ and $\mathbf{\Upsilon}$ all happen to be block-diagonal, which is the case in our problem, this projection matrix clearly is block-diagonal, with diagonal blocks $\mathrm{proj}((\mathbf{\Gamma}^{I}_{\boldsymbol{\vartheta};f_1})^{1/2}(\boldsymbol{\nu}^{(n)}_I)^{-1}\mathbf{\Upsilon}_I)$, $\mathrm{proj}((\mathbf{\Gamma}^{II}_{\boldsymbol{\vartheta};f_1})^{1/2}(\boldsymbol{\nu}^{(n)}_{II})^{-1}\mathbf{\Upsilon}_{II})$, and $\mathrm{proj}((\mathbf{\Gamma}^{III}_{\boldsymbol{\vartheta};f_1})^{1/2}(\boldsymbol{\nu}^{(n)}_{III})^{-1}\mathbf{\Upsilon}_{III})$, denoting projections in $\mathbb{R}^{mk}$, $\mathbb{R}^m$ and $\mathbb{R}^{mk_0}$, respectively. Since moreover $\mathcal{M}((\mathbf{\Gamma}^{I}_{\boldsymbol{\vartheta};f_1})^{1/2}(\boldsymbol{\nu}^{(n)}_I)^{-1} \times \mathbf{\Upsilon}_I) = \mathbb{R}^{mk}$, we have $\mathrm{proj}((\mathbf{\Gamma}^{I}_{\boldsymbol{\vartheta};f_1})^{1/2}(\boldsymbol{\nu}^{(n)}_I)^{-1}\mathbf{\Upsilon}_I) = \mathbf{I}_{mk}$, so that $Q_{\mathbf{\Upsilon}}$ reduces to

$$\begin{aligned}Q_{\mathbf{\Upsilon}} &= (\mathbf{\Delta}^{II}_{\boldsymbol{\vartheta};f_1})'(\mathbf{\Gamma}^{II}_{\boldsymbol{\vartheta};f_1})^{-1/2} \\ &\quad \times [\mathbf{I} - \mathrm{proj}((\mathbf{\Gamma}^{II}_{\boldsymbol{\vartheta};f_1})^{1/2}(\boldsymbol{\nu}^{(n)}_{II})^{-1}\mathbf{\Upsilon}_{II})](\mathbf{\Gamma}^{II}_{\boldsymbol{\vartheta};f_1})^{-1/2}\mathbf{\Delta}^{II}_{\boldsymbol{\vartheta};f_1} \\ &\quad + (\mathbf{\Delta}^{III}_{\boldsymbol{\vartheta};f_1})'(\mathbf{\Gamma}^{III}_{\boldsymbol{\vartheta};f_1})^{-1/2} \\ &\qquad \times [\mathbf{I} - \mathrm{proj}((\mathbf{\Gamma}^{III}_{\boldsymbol{\vartheta};f_1})^{1/2}(\boldsymbol{\nu}^{(n)}_{III})^{-1}\mathbf{\Upsilon}_{III})](\mathbf{\Gamma}^{III}_{\boldsymbol{\vartheta};f_1})^{-1/2}\mathbf{\Delta}^{III}_{\boldsymbol{\vartheta};f_1},\end{aligned} \tag{4.3}$$

where $\mathbf{\Delta}^{I}_{\boldsymbol{\vartheta};f_1}$ does not play any role. Accordingly, in the next section, we proceed with rank-based analogues of $\mathbf{\Delta}^{II}_{\boldsymbol{\vartheta};f_1}$ and $\mathbf{\Delta}^{III}_{\boldsymbol{\vartheta};f_1}$ only.

Note that the decomposition (4.3) of $Q_{\mathbf{\Upsilon}}$ into two asymptotically orthogonal quadratic forms corresponds to the decomposition of scatter heterogeneity into scale and shape heterogeneity; each of the two quadratic forms in the right-hand side of (4.3) actually constitutes a locally asymptotically optimal test statistic for $\mathcal{H}_0$ against one of those two subalternatives.

4.2. *A rank-based central sequence for scale and shape* (*scatter*). A general result by [20] implies that, in adaptive models for which fixed-$f_1$ submodels are ULAN and fixed-$\boldsymbol{\vartheta}$ submodels are generated by a group $\mathcal{G}_{\boldsymbol{\vartheta}}$, invariant versions of central sequences exist under very general assumptions. In the present context, this result would imply the existence, for the null values of $\boldsymbol{\vartheta}$, of central sequences based on the multivariate signs $\mathbf{U}_{ij}$ and the ranks $R_{ij}$. Although that result does not directly apply here, it is very likely that it still holds. This fact is confirmed by the asymptotic representation of Lemma 4.1(i) below.

Consider the signed-rank statistic [associated with some score function $K$ satisfying Assumption (C)] $\underset{\sim}{\boldsymbol{\Delta}}_{\boldsymbol{\vartheta};K} := ((\underset{\sim}{\boldsymbol{\Delta}}^{II}_{\boldsymbol{\vartheta};K})', (\underset{\sim}{\boldsymbol{\Delta}}^{III}_{\boldsymbol{\vartheta};K})')' := (\underset{\sim}{\Delta}^{II,1}_{\boldsymbol{\vartheta};K}, \ldots, \underset{\sim}{\Delta}^{II,m}_{\boldsymbol{\vartheta};K}, (\underset{\sim}{\boldsymbol{\Delta}}^{III,1}_{\boldsymbol{\vartheta};K})', \ldots, (\underset{\sim}{\boldsymbol{\Delta}}^{III,m}_{\boldsymbol{\vartheta};K})')'$, where

$$\underset{\sim}{\Delta}^{II,i}_{\boldsymbol{\vartheta};K} := \frac{n_i^{-1/2}}{2\sigma_i^2} \sum_{j=1}^{n_i} \left( K\left( \frac{R_{ij}}{n+1} \right) - k \right) \tag{4.4}$$

and

$$\underset{\sim}{\boldsymbol{\Delta}}^{III,i}_{\boldsymbol{\vartheta};K} := \frac{n_i^{-1/2}}{2} \mathbf{M}_k(\mathbf{V}_i)(\mathbf{V}_i^{\otimes 2})^{-1/2} \sum_{j=1}^{n_i} K\left( \frac{R_{ij}}{n+1} \right) \operatorname{vec}(\mathbf{U}_{ij}\mathbf{U}'_{ij}). \tag{4.5}$$

The following lemma provides (i) an asymptotic representation and (ii) the asymptotic distribution of $\underset{\sim}{\boldsymbol{\Delta}}_{\boldsymbol{\vartheta};K}$ (see the Appendix for the proof). An immediate corollary of (i) is that $\underset{\sim}{\boldsymbol{\Delta}}_{\boldsymbol{\vartheta};f_1} := \underset{\sim}{\boldsymbol{\Delta}}_{\boldsymbol{\vartheta};K_{f_1}}$, with $K = K_{f_1}$, actually constitutes a signed-rank version of the scatter part $((\boldsymbol{\Delta}^{II}_{\boldsymbol{\vartheta};f_1})', (\boldsymbol{\Delta}^{III}_{\boldsymbol{\vartheta};f_1})')'$ of the central sequence $\boldsymbol{\Delta}_{\boldsymbol{\vartheta};f_1}$.

LEMMA 4.1. *Assume that* (A), (B) *and* (C) *hold. Fix* $\boldsymbol{\vartheta} \in \mathcal{M}(\boldsymbol{\Upsilon})$ (*with common values* $\sigma$ *and* $\mathbf{V}$ *of the scale and shape parameters*). *Let* $R_{ij}$ *be the rank of* $d_{ij} := d_{ij}(\boldsymbol{\theta}_i, \mathbf{V})$ *among* $d_{11}, \ldots, d_{mn_m}$, *and let* $\mathbf{U}_{ij} := \mathbf{U}_{ij}(\boldsymbol{\theta}_i, \mathbf{V})$. *Then, for all* $g_1 \in \mathcal{F}_1$:

(i) $\underset{\sim}{\boldsymbol{\Delta}}_{\boldsymbol{\vartheta};K} = \boldsymbol{\Delta}_{\boldsymbol{\vartheta};K;g_1} + o_{L^2}(1)$, *under* $\mathrm{P}^{(n)}_{\boldsymbol{\vartheta};g_1}$, *as* $n \to \infty$, *where* $\boldsymbol{\Delta}_{\boldsymbol{\vartheta};K;g_1} := ((\boldsymbol{\Delta}^{II}_{\boldsymbol{\vartheta};K;g_1})', (\boldsymbol{\Delta}^{III}_{\boldsymbol{\vartheta};K;g_1})')' := (\Delta^{II,1}_{\boldsymbol{\vartheta};K;g_1}, \ldots, \Delta^{II,m}_{\boldsymbol{\vartheta};K;g_1}, (\boldsymbol{\Delta}^{III,1}_{\boldsymbol{\vartheta};K;g_1})', \ldots, (\boldsymbol{\Delta}^{III,m}_{\boldsymbol{\vartheta};K;g_1})')'$, *with*

$$\Delta^{II,i}_{\boldsymbol{\vartheta};K;g_1} := \frac{n_i^{-1/2}}{2\sigma^2} \sum_{j=1}^{n_i} \left( K\left( \tilde{G}_{1k}\left( \frac{d_{ij}}{\sigma} \right) \right) - k \right)$$

*and*

$$\text{(4.6)}\quad \mathbf{\Delta}^{III,i}_{\boldsymbol{\vartheta};K;g_1} := \frac{n_i^{-1/2}}{2}\mathbf{M}_k(\mathbf{V})(\mathbf{V}^{\otimes 2})^{-1/2}\sum_{j=1}^{n_i} K\left(\tilde{G}_{1k}\left(\frac{d_{ij}}{\sigma}\right)\right)\operatorname{vec}(\mathbf{U}_{ij}\mathbf{U}'_{ij});$$

(ii) *defining* $\mathbf{H}_k(\mathbf{V}) := \frac{1}{4k(k+2)}\mathbf{M}_k(\mathbf{V})[\mathbf{I}_{k^2} + \mathbf{K}_k](\mathbf{V}^{\otimes 2})^{-1}(\mathbf{M}_k(\mathbf{V}))'$, $\mathbf{\Delta}_{\boldsymbol{\vartheta};K;g_1}$ *is asymptotically normal with mean zero and mean*

$$\begin{pmatrix} \frac{1}{4\sigma^4}\mathcal{L}_k(K,g_1)\boldsymbol{\tau}_{II} \\ \mathcal{J}_k(K,g_1)[\mathbf{I}_m \otimes \mathbf{H}_k(\mathbf{V})]\boldsymbol{\tau}_{III} \end{pmatrix}$$

*under* $\mathrm{P}^{(n)}_{\boldsymbol{\vartheta};g_1}$ *and* $\mathrm{P}^{(n)}_{\boldsymbol{\vartheta}+n^{-1/2}\boldsymbol{\nu}^{(n)}\boldsymbol{\tau};g_1}$, *respectively, and covariance matrix*

$$\mathbf{\Gamma}_{\boldsymbol{\vartheta};K} := \operatorname{diag}(\mathbf{\Gamma}^{II}_{\boldsymbol{\vartheta};K}, \mathbf{\Gamma}^{III}_{\boldsymbol{\vartheta};K}) := \operatorname{diag}\left(\frac{1}{4\sigma^4}\mathcal{L}_k(K)\mathbf{I}_m, \mathcal{J}_k(K)[\mathbf{I}_m \otimes \mathbf{H}_k(\mathbf{V})]\right)$$

*under both (the claim under* $\mathrm{P}^{(n)}_{\boldsymbol{\vartheta}+n^{-1/2}\boldsymbol{\nu}^{(n)}\boldsymbol{\tau};g_1}$ *further requires* $g_1 \in \mathcal{F}_a$*).*

As mentioned in the description of the most stringent tests (see the comments after Proposition 4.1), we will need replacing the parameter $\boldsymbol{\vartheta}$ with some estimate. For this purpose, we assume the existence of $\hat{\boldsymbol{\vartheta}} := \hat{\boldsymbol{\vartheta}}^{(n)}$ satisfying:

Assumption (D). The sequence of estimators $(\hat{\boldsymbol{\vartheta}}^{(n)}, n \in \mathbb{N})$ is:

(D1) *constrained*: $\mathrm{P}^{(n)}_{\boldsymbol{\vartheta};g_1}[\hat{\boldsymbol{\vartheta}}^{(n)} \in \mathcal{M}(\mathbf{\Upsilon})] = 1$ for all $n$, $\boldsymbol{\vartheta} \in \mathcal{M}(\mathbf{\Upsilon})$, and $g_1 \in \mathcal{F}_1$;

(D2) $n^{1/2}(\boldsymbol{\nu}^{(n)})^{-1}$-*consistent*: for all $\boldsymbol{\vartheta} \in \mathcal{M}(\mathbf{\Upsilon})$, $n^{1/2}(\boldsymbol{\nu}^{(n)})^{-1}(\hat{\boldsymbol{\vartheta}}^{(n)} - \boldsymbol{\vartheta}) = O_{\mathrm{P}}(1)$, as $n \to \infty$, under $\bigcup_{g_1 \in \mathcal{F}_1}\{\mathrm{P}^{(n)}_{\boldsymbol{\vartheta};g_1}\}$;

(D3) *locally asymptotically discrete*: for all $\boldsymbol{\vartheta} \in \mathcal{M}(\mathbf{\Upsilon})$ and all $c > 0$, there exists $M = M(c) > 0$ such that the number of possible values of $\hat{\boldsymbol{\vartheta}}^{(n)}$ in balls of the form $\{\mathbf{t} \in \mathbb{R}^L : n^{1/2}\|(\boldsymbol{\nu}^{(n)})^{-1}(\mathbf{t} - \boldsymbol{\vartheta})\| \le c\}$ is bounded by $M$ as $n \to \infty$, and

(D4) *affine-equivariant*: denoting by $\hat{\boldsymbol{\vartheta}}^{(n)}(\mathbf{A}, \mathbf{B})$ the value of $\hat{\boldsymbol{\vartheta}}^{(n)}$ computed from the transformed sample $\mathbf{A}\mathbf{X}_{ij} + \mathbf{b}_i$, $j = 1, \ldots, n_i$, $i = 1, \ldots, m$, $\hat{\boldsymbol{\vartheta}}^{(n)}(\mathbf{A}, \mathbf{B}) = \mathbf{g}^{m,k}_{\mathbf{A},\mathbf{B}}(\hat{\boldsymbol{\vartheta}}^{(n)})$, for all $\mathbf{g}^{m,k}_{\mathbf{A},\mathbf{B}} \in \tilde{\mathcal{G}}^{m,k}$.

There are many possible choices for $\hat{\boldsymbol{\vartheta}}$. However, still in order to avoid moment assumptions, we propose the following estimators, related with the affine-equivariant median proposed by [23]. For each $i = 1, \ldots, m$, let $\hat{\boldsymbol{\theta}}_i$ and $\hat{\mathbf{V}}_i$ be characterized by

$$\frac{1}{n_i}\sum_{j=1}^{n_i}\mathbf{U}_{ij}(\hat{\boldsymbol{\theta}}_i, \hat{\mathbf{V}}_i) = \mathbf{0} \quad \text{and} \quad \frac{1}{n_i}\sum_{j=1}^{n_i}\mathbf{U}_{ij}(\hat{\boldsymbol{\theta}}_i, \hat{\mathbf{V}}_i)\mathbf{U}'_{ij}(\hat{\boldsymbol{\theta}}_i, \hat{\mathbf{V}}_i) = \frac{1}{k}\mathbf{I}_k,$$

with $|\hat{\mathbf{V}}_i| = 1$. Then, under $\boldsymbol{\vartheta} \in \mathcal{M}(\boldsymbol{\Upsilon})$, the common value $\mathbf{V}$ of the $\mathbf{V}_i$'s is consistently estimated [as $n \to \infty$, under $\bigcup_{g_1 \in \mathcal{F}_1} \{\mathrm{P}^{(n)}_{\boldsymbol{\vartheta};g_1}\}$ and Assumptions (A) and (B), *without any moment assumption on* $g_1$], at the rate required by Assumption (D2), by the Tyler [39] estimator $\hat{\mathbf{V}}$ computed from the $n$ data points $\mathbf{X}_{ij} - \hat{\boldsymbol{\theta}}_i$ and normalized in such a way that $|\hat{\mathbf{V}}| = 1$. Under the same conditions, the null value $\sigma$ of the scale is the common median of the i.i.d. radial distances $d_{ij}(\boldsymbol{\theta}_i, \mathbf{V})$, so that the empirical median $\hat{\sigma}$ of the $d_{ij}(\hat{\boldsymbol{\theta}}_i, \hat{\mathbf{V}})$'s can be used as an estimator of $\sigma$. Consequently, the estimator

$$\hat{\boldsymbol{\vartheta}} := (\hat{\boldsymbol{\theta}}'_1, \ldots, \hat{\boldsymbol{\theta}}'_m, \hat{\sigma}^2 \mathbf{1}'_m, \mathbf{1}'_m \otimes (\overset{\circ}{\mathrm{vech}}\, \hat{\mathbf{V}})')' \tag{4.7}$$

satisfies (D2) above—except perhaps for the $\hat{\boldsymbol{\vartheta}}_{II}$ part which, however, is not involved in the test statistics below. This estimator also satisfies (D1) and (D4). As for (D3), it is a purely technical requirement, with little practical implications (for fixed sample size, any estimator indeed can be considered part of a locally asymptotically discrete sequence). Therefore, we henceforth assume that (4.7) satisfies Assumption (D).

The resulting ranks $\hat{R}_{ij} := R_{ij}(\hat{\boldsymbol{\theta}}_1, \ldots, \hat{\boldsymbol{\theta}}_m, \hat{\mathbf{V}}, \ldots, \hat{\mathbf{V}})$ are usually called *aligned ranks*. The following *asymptotic linearity* result describes the asymptotic behavior of the aligned version $\underset{\sim}{\boldsymbol{\Delta}}_{\hat{\boldsymbol{\vartheta}};K}$ of $\underset{\sim}{\boldsymbol{\Delta}}_{\boldsymbol{\vartheta};K}$ under $\mathrm{P}^{(n)}_{\boldsymbol{\vartheta};g_1}$; see the Appendix for the proof.

PROPOSITION 4.2. *Assume that Assumptions* (A), (B), (C) *and* (D1)–(D3) *hold. Let* $g_1 \in \mathcal{F}_a$ *and* $\boldsymbol{\vartheta} \in \mathcal{M}(\boldsymbol{\Upsilon})$ *(with common values* $\sigma$ *and* $\mathbf{V}$ *for the scale and shape parameters). Then,*

$$\underset{\sim}{\boldsymbol{\Delta}}^{II}_{\hat{\boldsymbol{\vartheta}};K} - \underset{\sim}{\boldsymbol{\Delta}}^{II}_{\boldsymbol{\vartheta};K} + \frac{1}{4\sigma^4}\mathcal{L}_k(K, g_1)(\boldsymbol{\nu}^{(n)}_{II})^{-1} n^{1/2}(\hat{\boldsymbol{\vartheta}}^{(n)}_{II} - \boldsymbol{\vartheta}_{II})$$

*and*

$$\underset{\sim}{\boldsymbol{\Delta}}^{III}_{\hat{\boldsymbol{\vartheta}};K} - \underset{\sim}{\boldsymbol{\Delta}}^{III}_{\boldsymbol{\vartheta};K} + \mathcal{J}_k(K, g_1)[\mathbf{I}_m \otimes \mathbf{H}_k(\mathbf{V})](\boldsymbol{\nu}^{(n)}_{III})^{-1} n^{1/2}(\hat{\boldsymbol{\vartheta}}^{(n)}_{III} - \boldsymbol{\vartheta}_{III})$$

*are* $o_{\mathrm{P}}(1)$ *under* $\mathrm{P}^{(n)}_{\boldsymbol{\vartheta};g_1}$, *as* $n \to \infty$.

## 5. Optimal tests of scatter homogeneity.

5.1. *Optimal rank-based tests.* For all $\boldsymbol{\vartheta} \in \mathcal{M}(\boldsymbol{\Upsilon})$ (with common values $\sigma$ and $\mathbf{V}$ for the scale and shape parameters), define

$$\begin{aligned}\mathbf{P}^{II}_{\boldsymbol{\vartheta};K} &:= (\boldsymbol{\Gamma}^{II}_{\boldsymbol{\vartheta};K})^{-1} - (\boldsymbol{\nu}^{(n)}_{II})^{-1}\boldsymbol{\Upsilon}_{II}(\boldsymbol{\Upsilon}'_{II}(\boldsymbol{\nu}^{(n)}_{II})^{-1}\boldsymbol{\Gamma}^{II}_{\boldsymbol{\vartheta};K}(\boldsymbol{\nu}^{(n)}_{II})^{-1}\boldsymbol{\Upsilon}_{II})^{-1}\boldsymbol{\Upsilon}'_{II}(\boldsymbol{\nu}^{(n)}_{II})^{-1} \\ &= \frac{4\sigma^4}{\mathcal{L}_k(K)}[\mathbf{I}_m - \mathbf{C}^{(n)}]\end{aligned}$$

and

$$\begin{aligned}\mathbf{P}^{III}_{\boldsymbol{\vartheta};K} &:= (\boldsymbol{\Gamma}^{III}_{\boldsymbol{\vartheta};K})^{-1} - (\boldsymbol{\nu}^{(n)}_{III})^{-1}\boldsymbol{\Upsilon}_{III}(\boldsymbol{\Upsilon}'_{III}(\boldsymbol{\nu}^{(n)}_{III})^{-1}\boldsymbol{\Gamma}^{III}_{\boldsymbol{\vartheta};K}(\boldsymbol{\nu}^{(n)}_{III})^{-1}\boldsymbol{\Upsilon}_{III})^{-1}\boldsymbol{\Upsilon}'_{III}(\boldsymbol{\nu}^{(n)}_{III})^{-1}\\ &= (\mathcal{J}_k(K))^{-1}[\mathbf{I}_m - \mathbf{C}^{(n)}] \otimes (\mathbf{H}_k(\mathbf{V}))^{-1},\end{aligned}$$

where $\mathbf{C}^{(n)} = (C^{(n)}_{ii'})$ denotes the $(m\times m)$ matrix with entries $C^{(n)}_{ii'} := (\lambda^{(n)}_i\lambda^{(n)}_{i'})^{1/2}$. The $K$-score version $\underset{\sim}{\phi}{}^{(n)}_K$ of the rank-based tests we are proposing rejects $\mathcal{H}_0 : \boldsymbol{\vartheta} \in \mathcal{M}(\boldsymbol{\Upsilon})$ whenever the (affine-invariant) test statistic

$$\begin{aligned}\underset{\sim}{Q}{}^{(n)}_K &:= (\underset{\sim}{\boldsymbol{\Delta}}{}^{II}_{\hat{\boldsymbol{\vartheta}};K})'\mathbf{P}^{II}_{\hat{\boldsymbol{\vartheta}};K}\underset{\sim}{\boldsymbol{\Delta}}{}^{II}_{\hat{\boldsymbol{\vartheta}};K} + (\underset{\sim}{\boldsymbol{\Delta}}{}^{III}_{\hat{\boldsymbol{\vartheta}};K})'\mathbf{P}^{III}_{\hat{\boldsymbol{\vartheta}};K}\underset{\sim}{\boldsymbol{\Delta}}{}^{III}_{\hat{\boldsymbol{\vartheta}};K}\\ &= \sum_{i,i'=1}^{m}\left[\frac{\delta_{i,i'}}{n_i} - \frac{1}{n}\right]\\ &\quad\times \sum_{j=1}^{n_i}\sum_{j'=1}^{n_{i'}}\Biggl\{\frac{1}{\mathcal{L}_k(K)}\left(K\left(\frac{\hat{R}_{ij}}{n+1}\right) - k\right)\left(K\left(\frac{\hat{R}_{i'j'}}{n+1}\right) - k\right)\\ &\qquad + \frac{k(k+2)}{2\mathcal{J}_k(K)}K\left(\frac{\hat{R}_{ij}}{n+1}\right)K\left(\frac{\hat{R}_{i'j'}}{n+1}\right)\\ &\qquad\qquad \times\left((\hat{\mathbf{U}}'_{ij}\hat{\mathbf{U}}_{i'j'})^2 - \frac{1}{k}\right)\Biggr\}\end{aligned} \tag{5.1}$$

exceeds the $\alpha$-upper quantile $\chi^2_{(m-1)(k_0+1);1-\alpha}$ of the chi-square distribution with $(m-1)(k_0+1)$ degrees of freedom ($\delta_{i,i'}$ stands for the usual Kronecker symbol); the explicit form of $(\mathbf{H}_k(\mathbf{V}))^{-1}$ allowing for (5.1) can be found in Lemma 5.2 of [18]. In the sequel, we write $\underset{\sim}{\phi}{}^{(n)}_{f_1}$ and $\underset{\sim}{Q}{}^{(n)}_{f_1}$ for $\underset{\sim}{\phi}{}^{(n)}_{K_{f_1}}$ and $\underset{\sim}{Q}{}^{(n)}_{K_{f_1}}$, respectively.

The decomposition (5.1) of the rank-based test statistic $\underset{\sim}{Q}{}^{(n)}_K$ into two asymptotically orthogonal terms parallels the corresponding decomposition (4.3) of $\mathbf{Q}_{\boldsymbol{\Upsilon}}$ (see the closing remark of Section 4.1), with the same interpretation in terms of subalternatives of scale and shape heterogeneity, respectively.

We are now ready to state the main result of this paper; for the sake of simplicity, asymptotic powers are expressed under Assumption (B′) and perturbations $\boldsymbol{\tau}^{(n)}$ such that $\lim_{n\to\infty}\boldsymbol{\nu}^{(n)}\boldsymbol{\tau}^{(n)} = \boldsymbol{\nu\tau} \notin \mathcal{M}(\boldsymbol{\Upsilon})$, with $\boldsymbol{\nu}_{II}\boldsymbol{\tau}_{II} = (s_1^2/\sqrt{\lambda_1}, \ldots, s_m^2/\sqrt{\lambda_m})'$ and $\boldsymbol{\nu}_{III}\boldsymbol{\tau}_{III} = ((\overset{\circ}{\mathrm{vech}}\,\mathbf{v}_1)'/\sqrt{\lambda_1}, \ldots, (\overset{\circ}{\mathrm{vech}}\,\mathbf{v}_m)'/\sqrt{\lambda_m})'$. For any such $\boldsymbol{\tau}$ and any $\boldsymbol{\vartheta} \in \mathcal{M}(\boldsymbol{\Upsilon})$ (still with common values $\sigma$ and $\mathbf{V}$ of the scale and shape parameters), let

$$r^{II}_{\boldsymbol{\vartheta},\boldsymbol{\tau}} := \frac{1}{\sigma^4}\lim_{n\to\infty}\{(\boldsymbol{\tau}^{(n)}_{II})'[\mathbf{I}_m - \mathbf{C}^{(n)}]\boldsymbol{\tau}^{(n)}_{II}\}$$

$$= \sum_{1\le i<i'\le m} \frac{\lambda_i\lambda_{i'}}{\sigma^4}\left(\frac{s_i^2}{\sqrt{\lambda_i}} - \frac{s_{i'}^2}{\sqrt{\lambda_{i'}}}\right)^2 \tag{5.2}$$

and

$$\begin{aligned} r^{III}_{\boldsymbol{\vartheta},\boldsymbol{\tau}} &:= 2k(k+2)\lim_{n\to\infty}\{(\boldsymbol{\tau}^{(n)}_{III})'[[\mathbf{I}_m - \mathbf{C}^{(n)}]\otimes \mathbf{H}_k(\mathbf{V})]\boldsymbol{\tau}^{(n)}_{III}\} \\ &= \sum_{1\le i<i'\le m} \lambda_i\lambda_{i'} \operatorname{tr}\left[\left(\mathbf{V}^{-1}\left(\frac{\mathbf{v}_i}{\sqrt{\lambda_i}} - \frac{\mathbf{v}_{i'}}{\sqrt{\lambda_{i'}}}\right)\right)^2\right]; \end{aligned} \tag{5.3}$$

note that $\operatorname{tr}(\mathbf{V}^{-1}\mathbf{v}_i) = 0$ for all $i$ (see the comments before Proposition 4.1).

THEOREM 5.1. *Assume that* (A), (B), (C) *and* (D1)–(D3) *hold. Then:*

(i) $\underset{\sim}{Q}^{(n)}_K$ *is asymptotically chi-square with* $(m-1)(k_0+1)$ *degrees of freedom under* $\bigcup_{\boldsymbol{\vartheta}\in\mathcal{M}(\boldsymbol{\Upsilon})}\bigcup_{g_1\in\mathcal{F}_a}\{\mathrm{P}^{(n)}_{\boldsymbol{\vartheta};g_1}\}$*, and [provided that* (B) *is reinforced into* (B′)*] asymptotically noncentral chi-square, still with* $(m-1)(k_0+1)$ *degrees of freedom, but with noncentrality parameter*

$$\frac{\mathcal{L}^2_k(K,g_1)}{4\mathcal{L}_k(K)} r^{II}_{\boldsymbol{\vartheta},\boldsymbol{\tau}} + \frac{\mathcal{J}^2_k(K,g_1)}{2k(k+2)\mathcal{J}_k(K)} r^{III}_{\boldsymbol{\vartheta},\boldsymbol{\tau}} \tag{5.4}$$

*under* $\mathrm{P}^{(n)}_{\boldsymbol{\vartheta}+n^{-1/2}\boldsymbol{\nu}^{(n)}\boldsymbol{\tau}^{(n)};g_1}$, $\boldsymbol{\vartheta}\in\mathcal{M}(\boldsymbol{\Upsilon})$, $\boldsymbol{\nu\tau} := \lim_{n\to\infty}\boldsymbol{\nu}^{(n)}\boldsymbol{\tau}^{(n)} \notin \mathcal{M}(\boldsymbol{\Upsilon})$*, and* $g_1\in\mathcal{F}_a$*;*

(ii) *the sequence of tests* $\underset{\sim}{\phi}^{(n)}_K$ *has asymptotic level* $\alpha$ *under* $\bigcup_{\boldsymbol{\vartheta}\in\mathcal{M}(\boldsymbol{\Upsilon})}\bigcup_{g_1\in\mathcal{F}_a}\{\mathrm{P}^{(n)}_{\boldsymbol{\vartheta};g_1}\}$*;*

(iii) *if* $f_1\in\mathcal{F}_a$ *and* $K_{f_1}$ *satisfies Assumption* (C)*, the sequence of tests* $\underset{\sim}{\phi}^{(n)}_{f_1}$ *is locally and asymptotically most stringent, still at asymptotic level* $\alpha$*, for* $\bigcup_{\boldsymbol{\vartheta}\in\mathcal{M}(\boldsymbol{\Upsilon})}\bigcup_{g_1\in\mathcal{F}_a}\{\mathrm{P}^{(n)}_{\boldsymbol{\vartheta};g_1}\}$ *against alternatives of the form* $\bigcup_{\boldsymbol{\vartheta}\notin\mathcal{M}(\boldsymbol{\Upsilon})}\{\mathrm{P}^{(n)}_{\boldsymbol{\vartheta};f_1}\}$.

See the Appendix for the proof. After some algebra, one obtains

$$\underset{\sim}{Q}^{(n)}_K = \frac{1}{n}\sum_{1\le i<i'\le m}(n_i+n_{i'})\underset{\sim}{Q}^{(n)}_{K;i,i'},$$

where

$$\underset{\sim}{Q}^{(n)}_{K;i,i'} = \frac{n_i n_{i'}}{n_i+n_{i'}}$$

$$
\begin{aligned}
&\times \left\{ \frac{1}{\mathcal{L}_k(K)} \left[ \frac{1}{n_i} \sum_{j=1}^{n_i} K\left(\frac{\hat{R}_{ij}}{n+1}\right) - \frac{1}{n_{i'}} \sum_{j'=1}^{n_{i'}} K\left(\frac{\hat{R}_{i'j'}}{n+1}\right) \right]^2 \right. \\
&\quad + \frac{k(k+2)}{2\mathcal{J}_k(K)} \\
&\quad \times \left\| \left[ \frac{1}{n_i} \sum_{j=1}^{n_i} K\left(\frac{\hat{R}_{ij}}{n+1}\right) \operatorname{vec}\left(\hat{\mathbf{U}}_{ij}\hat{\mathbf{U}}'_{ij} - \frac{1}{k}\mathbf{I}_k\right) \right] \right. \\
&\quad \left.\left. - \left[ \frac{1}{n_{i'}} \sum_{j'=1}^{n_{i'}} K\left(\frac{\hat{R}_{i'j'}}{n+1}\right) \operatorname{vec}\left(\hat{\mathbf{U}}_{i'j'}\hat{\mathbf{U}}'_{i'j'} - \frac{1}{k}\mathbf{I}_k\right) \right] \right\|^2 \right\}
\end{aligned}
\tag{5.5}
$$

is the test statistic obtained in the two-sample case (for populations $i$ and $i'$); see [40] for a similar decomposition in MANOVA problems. As announced, no estimate $\hat{\boldsymbol{\vartheta}}_{II}$ of the common scale appears in the test statistics. Also, letting

$$
\underset{\sim}{\mathbf{S}}_{K;i} := \frac{1}{n_i} \sum_{j=1}^{n_i} K\left(\frac{\hat{R}_{ij}}{n+1}\right) \hat{\mathbf{U}}_{ij}\hat{\mathbf{U}}'_{ij},
$$

the statistics in (5.5) take the simple form

$$
\begin{aligned}
\underset{\sim}{Q}^{(n)}_{K;i,i'} &= \frac{n_i n_{i'}}{n_i + n_{i'}} \left\{ \frac{1}{\mathcal{L}_k(K)} \operatorname{tr}^2[\underset{\sim}{\mathbf{S}}_{K;i} - \underset{\sim}{\mathbf{S}}_{K;i'}] \right. \\
&\quad \left. + \frac{k(k+2)}{2\mathcal{J}_k(K)} \left[ \operatorname{tr}[(\underset{\sim}{\mathbf{S}}_{K;i} - \underset{\sim}{\mathbf{S}}_{K;i'})^2] - \frac{1}{k} \operatorname{tr}^2[\underset{\sim}{\mathbf{S}}_{K;i} - \underset{\sim}{\mathbf{S}}_{K;i'}] \right] \right\} \\
&= \frac{n_i n_{i'}}{n_i + n_{i'}} \left\{ \frac{k(k+2)}{2\mathcal{J}_k(K)} \operatorname{tr}[(\underset{\sim}{\mathbf{S}}_{K;i} - \underset{\sim}{\mathbf{S}}_{K;i'})^2] \right. \\
&\quad \left. - \frac{k(\mathcal{J}_k(K) - k(k+2))}{2\mathcal{J}_k(K)\mathcal{L}_k(K)} \operatorname{tr}^2[\underset{\sim}{\mathbf{S}}_{K;i} - \underset{\sim}{\mathbf{S}}_{K;i'}] \right\}.
\end{aligned}
$$

For Gaussian scores (i.e., for $K = K_{\phi_1}$; see Section 2.3), one obtains the *van der Waerden* test statistics

$$
\underset{\sim}{Q}^{(n)}_{\mathrm{vdW}} = \frac{1}{n} \sum_{1 \le i < i' \le m} (n_i + n_{i'}) \underset{\sim}{Q}^{(n)}_{\mathrm{vdW};i,i'},
\tag{5.6}
$$

where

$$
\underset{\sim}{Q}^{(n)}_{\mathrm{vdW};i,i'} = \frac{n_i n_{i'}}{2(n_i + n_{i'})} \operatorname{tr}[(\underset{\sim}{\mathbf{S}}_{\mathrm{vdW};i} - \underset{\sim}{\mathbf{S}}_{\mathrm{vdW};i'})^2],
$$

with $\underset{\sim}{\mathbf{S}}_{\mathrm{vdW};i} := n_i^{-1}\sum_{j=1}^{n_i}\Psi_k^{-1}(\hat{R}_{ij}/(n+1))\hat{\mathbf{U}}_{ij}\hat{\mathbf{U}}'_{ij}$. The Student scores (i.e., $K = K_{f^t_{1,\nu}}$; see Section 2.3 again) yield

$$\underset{\sim}{Q}^{(n)}_{f^t_{1,\nu}} = \frac{1}{n}\sum_{1\le i<i'\le m}(n_i+n_{i'})\underset{\sim}{Q}^{(n)}_{f^t_{1,\nu};i,i'}, \tag{5.7}$$

where

$$\underset{\sim}{Q}^{(n)}_{f^t_{1,\nu};i,i'} = \frac{n_i n_{i'}}{n_i+n_{i'}}\frac{k+\nu+2}{2(k+\nu)} \times\left\{\mathrm{tr}[(\underset{\sim}{\mathbf{S}}_{f^t_{1,\nu};i}-\underset{\sim}{\mathbf{S}}_{f^t_{1,\nu};i'})^2]+\frac{1}{\nu}\mathrm{tr}^2[\underset{\sim}{\mathbf{S}}_{f^t_{1,\nu};i}-\underset{\sim}{\mathbf{S}}_{f^t_{1,\nu};i'}]\right\}$$

with $\underset{\sim}{\mathbf{S}}_{f^t_{1,\nu};i} := k(k+\nu)n_i^{-1}\sum_{j=1}^{n_i}G^{-1}_{k,\nu}(\hat{R}_{ij}/(n+1))/[\nu+kG^{-1}_{k,\nu}(\hat{R}_{ij}/(n+1))]\times \hat{\mathbf{U}}_{ij}\hat{\mathbf{U}}'_{ij}$. As for the tests associated with the power score functions $K_a$, they are based on

$$\underset{\sim}{Q}^{(n)}_{K_a} = \frac{1}{n}\sum_{1\le i<i'\le m}(n_i+n_{i'})\underset{\sim}{Q}^{(n)}_{K_a;i,i'}, \qquad a>0, \tag{5.8}$$

where, letting $\underset{\sim}{\mathbf{S}}_{K_a;i} := k(a+1)(n+1)^{-a}n_i^{-1}\sum_{j=1}^{n_i}(\hat{R}_{ij})^a\hat{\mathbf{U}}_{ij}\hat{\mathbf{U}}'_{ij}$,

$$\underset{\sim}{Q}^{(n)}_{K_a;i,i'} = \frac{n_i n_{i'}}{n_i+n_{i'}}\frac{2a+1}{2a^2(a+1)^2k^2}\{a^2k(k+2)\,\mathrm{tr}[(\underset{\sim}{\mathbf{S}}_{K_a;i}-\underset{\sim}{\mathbf{S}}_{K_a;i'})^2] - (a^2k-4a-2)\,\mathrm{tr}^2[\underset{\sim}{\mathbf{S}}_{K_a;i}-\underset{\sim}{\mathbf{S}}_{K_a;i'}]\}.$$

COROLLARY 5.1. *Assume that the conditions of Theorem 5.1 hold. Then:*

(i) *provided that $g_1\in\mathcal{F}_a$ is such that $\mathcal{L}_k(K,g_1)\neq 0\neq\mathcal{J}_k(K,g_1)$, $\underset{\sim}{\phi}^{(n)}_K$ is consistent under any local $g_1$-alternative [i.e., under any $\mathrm{P}^{(n)}_{\boldsymbol{\vartheta}+n^{-1/2}\boldsymbol{\nu}^{(n)}\boldsymbol{\tau}^{(n)};g_1}$, $\boldsymbol{\vartheta}\in\mathcal{M}(\boldsymbol{\Upsilon})$, $\lim_{n\to\infty}\boldsymbol{\nu}^{(n)}\boldsymbol{\tau}^{(n)}\notin\mathcal{M}(\boldsymbol{\Upsilon})$];*

(ii) *the same conclusion holds if $u\mapsto K(u)$ absolutely continuous with a.e. derivative $\dot{K}$, and if $g_1\in\mathcal{F}_a$ is such that $\int_0^\infty\dot{K}(\tilde{G}_{1k}(r))r(\tilde{g}_{1k}(r))^2\,dr>0$ (in particular, if $\dot{K}$ is nondecreasing).*

We refer to the Appendix for the proof. This corollary shows that the van der Waerden tests above, as well as those achieving local asymptotic stringency at prespecified Student or power-exponential densities, are locally consistent against arbitrary elliptical alternatives, since the corresponding score functions are strictly increasing. Similar conclusions hold for the tests associated with the power score functions $K_a$, $a > 0$. More general consistency results, against nonlocal and possibly nonelliptical alternatives, of course, are highly desirable, and can be obtained by exploiting Hájek's results on rank statistics under alternatives (see [13]), much along the same lines as in Section 5.2 of [17]. Such results, which we do not include here, imply that consistency is achieved at "almost all" nonlocal alternatives—the exception being those very particular densities achieving a set of orthogonality conditions involving the scores and the $\mathbf{U}_{ij}$'s.

Another important concern is validity under homogeneous scatter but possibly heterogeneous or/and nonelliptical densities. The ranks of the $d_{ij}$'s then lose their distribution-freeness and/or their independence with respect to the $\mathbf{U}_{ij}$'s. Under totally arbitrary situations, little can be said, but it is extremely unlikely that the tests we are developing here remain valid; if ellipticity and finite moments of order four can be assumed, the pseudo-Gaussian methods developed in [19], which remain valid under possibly heterokurtic elliptic densities, are preferable.

5.2. *The optimal pseudo-Gaussian tests.* As explained in the Introduction, the pseudo-Gaussian tests $\phi^{(n)}_{\mathcal{N}*}$ and $\phi^{(n)}_{\mathrm{Schott}*}$ are the natural competitors of our rank-based procedures. Since they are asymptotically equivalent (see [19]), we concentrate on $\phi^{(n)}_{\mathcal{N}*}$. This test is valid under finite fourth-order moments $E_k(g_1) := \int_0^1 (\tilde{G}_{1k}^{-1}(u))^4\, du = (\mu_{k-1;g_1})^{-1} \int_0^\infty r^{k+3} g_1(r)\, dr$, that is, under any $\mathrm{P}^{(n)}_{\boldsymbol{\vartheta};g_1}$ such that $g_1 \in \mathcal{F}_1^4 := \{g_1 \in \mathcal{F}_1 : E_k(g_1) < \infty\}$. For all such $g_1$, let $D_k(g_1) := \int_0^1 (\tilde{G}_{1k}^{-1}(u))^2\, du$: then $\kappa_k(g_1) := \frac{k}{(k+2)} \frac{E_k(g_1)}{D_k^2(g_1)} - 1$ is a measure of kurtosis (see, e.g., page 54 of [1]), common to the $m$ elliptic populations under $\mathrm{P}^{(n)}_{\boldsymbol{\vartheta};g_1}$. This kurtosis is consistently (still under $\mathrm{P}^{(n)}_{\boldsymbol{\vartheta};g_1}$, $g_1 \in \mathcal{F}_1^4$) estimated by

$$\hat{\kappa}_k := \frac{k}{(k+2)} \frac{(n^{-1} \sum_{i=1}^m \sum_{j=1}^{n_i} \hat{d}_{ij}^4)}{(n^{-1} \sum_{i=1}^m \sum_{j=1}^{n_i} \hat{d}_{ij}^2)^2} - 1,$$

where $\hat{d}_{ij} := d_{ij}(\bar{\mathbf{X}}_i, \mathbf{S}_i)$, with $\bar{\mathbf{X}}_i := n_i^{-1} \sum_{j=1}^{n_i} \mathbf{X}_{ij}$ and $\mathbf{S}_i := n_i^{-1} \sum_{j=1}^{n_i} (\mathbf{X}_{ij} - \bar{\mathbf{X}}_i)(\mathbf{X}_{ij} - \bar{\mathbf{X}}_i)'$. At the multinormal ($g_1 = \phi_1$), $E_k(\phi_1) = k(k+2)/a_k^2$ and $D_k(\phi_1) = k/a_k$, so that $\kappa_k(\phi_1) = 0$.

The pseudo-Gaussian test $\phi^{(n)}_{\mathcal{N}*}$ rejects $\mathcal{H}_0$ (at asymptotic level $\alpha$) whenever

$$Q^{(n)}_{\mathcal{N}*} := \frac{1}{n} \sum_{1 \le i < i' \le m} (n_i + n_{i'}) Q^{(n)}_{\mathcal{N}*;i,i'} > \chi^2_{(m-1)(k_0+1);1-\alpha}, \tag{5.9}$$

where, letting $\mathbf{S} := n^{-1}\sum_{i=1}^m n_i \mathbf{S}_i$,

$$Q^{(n)}_{\mathcal{N}*;i,i'} := \frac{n_i n_{i'}}{n_i + n_{i'}} \frac{1}{2(1+\hat{\kappa}_k)} \times \left\{ \mathrm{tr}[(\mathbf{S}^{-1}(\mathbf{S}_i - \mathbf{S}_{i'}))^2] - \frac{\hat{\kappa}_k}{(k+2)\hat{\kappa}_k + 2} \mathrm{tr}^2[\mathbf{S}^{-1}(\mathbf{S}_i - \mathbf{S}_{i'})] \right\}.$$

This test statistic is clearly affine-invariant; the following Theorem (see [19]) summarizes its asymptotic properties, which also are those of $\phi^{(n)}_{\mathrm{Schott}*}$.

THEOREM 5.2. *Assume that (A) and (B) hold. Then:*

(i) $Q^{(n)}_{\mathcal{N}*}$ *is asymptotically chi-square with* $(m-1)(k_0+1)$ *degrees of freedom under* $\bigcup_{\boldsymbol{\vartheta}\in\mathcal{M}(\boldsymbol{\Upsilon})} \bigcup_{g_1\in\mathcal{F}_1^4} \{\mathrm{P}^{(n)}_{\boldsymbol{\vartheta};g_1}\}$, *and [provided that (B) is reinforced into (B′)] asymptotically noncentral chi-square, still with* $(m-1)(k_0+1)$ *degrees of freedom but with noncentrality parameter*

$$\frac{k}{(k+2)\kappa_k(g_1)+2} r^{II}_{\boldsymbol{\vartheta},\boldsymbol{\tau}} + \frac{1}{2(1+\kappa_k(g_1))} r^{III}_{\boldsymbol{\vartheta},\boldsymbol{\tau}} \tag{5.10}$$

*under* $\mathrm{P}^{(n)}_{\boldsymbol{\vartheta}+n^{-1/2}\boldsymbol{\nu}^{(n)}\boldsymbol{\tau}^{(n)};g_1}$, *with* $\boldsymbol{\vartheta}\in\mathcal{M}(\boldsymbol{\Upsilon})$, $\boldsymbol{\nu\tau} := \lim_{n\to\infty}\boldsymbol{\nu}^{(n)}\boldsymbol{\tau}^{(n)} \notin \mathcal{M}(\boldsymbol{\Upsilon})$, $g_1 \in \mathcal{F}_a^4 (:= \mathcal{F}_1^4 \cap \mathcal{F}_a)$, $r^{II}_{\boldsymbol{\vartheta},\boldsymbol{\tau}}$ *and* $r^{III}_{\boldsymbol{\vartheta},\boldsymbol{\tau}}$ *defined in (5.2) and (5.3);*

(ii) $\phi^{(n)}_{\mathcal{N}*}$ *has asymptotic level* $\alpha$ *under* $\bigcup_{\boldsymbol{\vartheta}\in\mathcal{M}(\boldsymbol{\Upsilon})} \bigcup_{g_1\in\mathcal{F}_1^4} \{\mathrm{P}^{(n)}_{\boldsymbol{\vartheta};g_1}\}$;

(iii) $\phi^{(n)}_{\mathcal{N}*}$ *is locally and asymptotically most stringent, still at asymptotic level* $\alpha$, *for* $\bigcup_{\boldsymbol{\vartheta}\in\mathcal{M}(\boldsymbol{\Upsilon})} \bigcup_{g_1\in\mathcal{F}_1^4} \{\mathrm{P}^{(n)}_{\boldsymbol{\vartheta};g_1}\}$ *against alternatives of the form* $\bigcup_{\boldsymbol{\vartheta}\notin\mathcal{M}(\boldsymbol{\Upsilon})} \{\mathrm{P}^{(n)}_{\boldsymbol{\vartheta};\phi_1}\}$.

**6. Asymptotic relative efficiencies.** The asymptotic relative efficiencies of the rank-based tests $\underset{\sim}{\phi}{}^{(n)}_K$ with respect to $\phi^{(n)}_{\mathcal{N}*}$ and $\phi^{(n)}_{\mathrm{Schott}*}$ directly follow as ratios of noncentrality parameters under local alternatives (see Theorems 5.1 and 5.2).

PROPOSITION 6.1. *Assume that (A), (B′), (C) and (D) hold, and that* $g_1 \in \mathcal{F}_a^4$. *Then, the asymptotic relative efficiency of* $\underset{\sim}{\phi}{}^{(n)}_K$ *with respect to* $\phi^{(n)}_{\mathcal{N}*}$, *when testing* $\mathrm{P}^{(n)}_{\boldsymbol{\vartheta};g_1}$ *against* $\mathrm{P}^{(n)}_{\boldsymbol{\vartheta}+n^{-1/2}\boldsymbol{\nu}^{(n)}\boldsymbol{\tau}^{(n)};g_1}$ *[*$\boldsymbol{\vartheta}\in\mathcal{M}(\boldsymbol{\Upsilon})$ *and* $\boldsymbol{\nu\tau} := \lim_{n\to\infty}\boldsymbol{\nu}^{(n)}\boldsymbol{\tau}^{(n)} \notin \mathcal{M}(\boldsymbol{\Upsilon})$*], is*

$$\mathrm{ARE}_{\boldsymbol{\vartheta},\boldsymbol{\tau},k,g_1}(\underset{\sim}{\phi}{}^{(n)}_K / \phi^{(n)}_{\mathcal{N}*}) = (1-\xi)\mathrm{ARE}^{(\mathrm{scale})}_{k,g_1}(\underset{\sim}{\phi}{}^{(n)}_K / \phi^{(n)}_{\mathcal{N}*}) + \xi \mathrm{ARE}^{(\mathrm{shape})}_{k,g_1}(\underset{\sim}{\phi}{}^{(n)}_K / \phi^{(n)}_{\mathcal{N}*}),$$

*where*

$$\text{ARE}^{(\text{scale})}_{k,g_1}(\underset{\sim}{\phi}^{(n)}_K/\phi^{(n)}_{\mathcal{N}*}) := \frac{((k+2)\kappa_k(g_1)+2)\mathcal{L}^2_k(K,g_1)}{4k\mathcal{L}_k(K)}, \tag{6.1}$$

$$\text{ARE}^{(\text{shape})}_{k,g_1}(\underset{\sim}{\phi}^{(n)}_K/\phi^{(n)}_{\mathcal{N}*}) := \frac{(1+\kappa_k(g_1))\mathcal{J}^2_k(K,g_1)}{k(k+2)\mathcal{J}_k(K)} \tag{6.2}$$

*and* $\xi := \xi_{\boldsymbol{\vartheta},\boldsymbol{\tau},k,g_1} \in [0,1]$ *is given by*

$$\begin{aligned}\xi_{\boldsymbol{\vartheta},\boldsymbol{\tau},k,g_1} &:= ((k+2)\kappa_k(g_1)+2)r^{III}_{\boldsymbol{\vartheta},\boldsymbol{\tau}}\\ &\quad\times[2k(1+\kappa_k(g_1))r^{II}_{\boldsymbol{\vartheta},\boldsymbol{\tau}}+((k+2)\kappa_k(g_1)+2)r^{III}_{\boldsymbol{\vartheta},\boldsymbol{\tau}}]^{-1}.\end{aligned}$$

The "shape AREs" in (6.2) coincide with those obtained in problems involving shape only—such as testing null hypotheses of the form $\mathbf{V}=\mathbf{V}_0$ for fixed $\mathbf{V}_0$ (see [17]). Proposition 6.1 shows that the AREs, with respect to the pseudo-Gaussian tests of Section 5.2, of the rank tests proposed in Section 5.1 are convex linear combinations of these "shape AREs" and the "scale AREs" in (6.1).

Numerical values of (6.1) and (6.2), for various values of the space dimension $k$ and various radial densities (Student, Gaussian and power-exponential), are given in Table 1 for the van der Waerden test $\underset{\sim}{\phi}^{(n)}_{\text{vdW}}$, the Wilcoxon test $\underset{\sim}{\phi}^{(n)}_{K_1}$, and the Spearman test $\underset{\sim}{\phi}^{(n)}_{K_2}$ (the score functions $K_a$, $a>0$ were defined in Section 2.3). These ARE values are uniformly large (with the exception, possibly, of univariate scale Wilcoxon AREs), particularly so under heavy tails, as often in rank-based inference. Also note that the AREs of the proposed van der Waerden tests with respect to the parametric Gaussian tests are larger than or equal to one for all distributions considered in Table 1. For pure shape alternatives, [29] has shown that a Chernoff–Savage property holds, that is, $\inf_{g_1}\text{ARE}^{(\text{shape})}_{k,g_1}(\underset{\sim}{\phi}^{(n)}_{\text{vdW}}/\phi^{(n)}_{\mathcal{N}*})=1$. One may wonder whether this uniform dominance property of van der Waerden tests extends to the present situation. Although it does for usual distributions, including all Student and power-exponential ones, the general answer unfortunately is negative; see Section 4 of [29] for a (pathological) counter example.

Note that Theorem 5.2 clearly shows that $\phi^{(n)}_{\mathcal{N}*}$ and $\phi^{(n)}_{\text{Schott}*}$ are not efficiency-robust. Indeed, the noncentrality parameter (5.10) under Student radial densities with $4+\delta$ degrees of freedom tends to zero as $\delta\to 0$, so that asymptotic local powers are arbitrarily close to the nominal level $\alpha$. The efficiency-robustness of our rank tests, quite on the contrary, is not affected, as the ARE values (6.1) and (6.2), under the same Student densities with $4+\delta$ degrees of freedom, both tend to infinity as $\delta\to 0$.

Table 1

*AREs, for $\xi = 0$ (pure scale alternatives) and $\xi = 1$ (pure shape alternatives), of the van der Waerden (vdW), Wilcoxon (W), and Spearman (SP) rank-based tests with respect to the pseudo- Gaussian tests, under k-dimensional Student (with 5, 8 and 12 degrees of freedom), Gaussian, and power-exponential densities (with parameter $\eta = 2, 3, 5$), for $k = 2, 3, 4, 6, 10$ and $k \to \infty$*

| | | | Underlying density | | | | | | |
|---|---|---|---|---|---|---|---|---|---|
| | $k$ | $\xi$ | $t_5$ | $t_8$ | $t_{12}$ | $\mathcal{N}$ | $e_2$ | $e_3$ | $e_5$ |
| vdW | 1 | 0 | 2.321 | 1.230 | 1.082 | 1.000 | 1.151 | 1.376 | 1.822 |
| | | 1 | —— | —— | —— | —— | —— | —— | —— |
| | 2 | 0 | 2.551 | 1.280 | 1.102 | 1.000 | 1.115 | 1.296 | 1.669 |
| | | 1 | 2.204 | 1.215 | 1.078 | 1.000 | 1.129 | 1.308 | 1.637 |
| | 3 | 0 | 2.732 | 1.322 | 1.120 | 1.000 | 1.092 | 1.241 | 1.558 |
| | | 1 | 2.270 | 1.233 | 1.086 | 1.000 | 1.108 | 1.259 | 1.536 |
| | 4 | 0 | 2.881 | 1.358 | 1.136 | 1.000 | 1.077 | 1.202 | 1.475 |
| | | 1 | 2.326 | 1.249 | 1.093 | 1.000 | 1.093 | 1.223 | 1.462 |
| | 6 | 0 | 3.108 | 1.416 | 1.163 | 1.000 | 1.057 | 1.151 | 1.361 |
| | | 1 | 2.413 | 1.275 | 1.106 | 1.000 | 1.072 | 1.174 | 1.363 |
| | 10 | 0 | 3.403 | 1.498 | 1.204 | 1.000 | 1.037 | 1.099 | 1.239 |
| | | 1 | 2.531 | 1.312 | 1.126 | 1.000 | 1.050 | 1.121 | 1.254 |
| | $\infty$ | 0 | 4.586 | 1.894 | 1.446 | 1.000 | 1.000 | 1.000 | 1.000 |
| | | 1 | 3.000 | 1.500 | 1.250 | 1.000 | 1.000 | 1.000 | 1.000 |
| W | 1 | 0 | 1.993 | 0.939 | 0.769 | 0.608 | 0.519 | 0.509 | 0.517 |
| | | 1 | —— | —— | —— | —— | —— | —— | —— |
| | 2 | 0 | 2.604 | 1.185 | 0.959 | 0.750 | 0.694 | 0.703 | 0.743 |
| | | 1 | 2.258 | 1.174 | 1.001 | 0.844 | 0.789 | 0.804 | 0.842 |
| | 3 | 0 | 2.929 | 1.304 | 1.045 | 0.811 | 0.775 | 0.795 | 0.854 |
| | | 1 | 2.386 | 1.246 | 1.068 | 0.913 | 0.897 | 0.933 | 1.001 |
| | 4 | 0 | 3.140 | 1.377 | 1.096 | 0.844 | 0.820 | 0.844 | 0.911 |
| | | 1 | 2.432 | 1.273 | 1.094 | 0.945 | 0.955 | 1.006 | 1.095 |
| | 6 | 0 | 3.407 | 1.467 | 1.156 | 0.879 | 0.866 | 0.892 | 0.961 |
| | | 1 | 2.451 | 1.283 | 1.105 | 0.969 | 1.008 | 1.075 | 1.188 |
| | 10 | 0 | 3.685 | 1.560 | 1.216 | 0.908 | 0.903 | 0.925 | 0.984 |
| | | 1 | 2.426 | 1.264 | 1.088 | 0.970 | 1.032 | 1.106 | 1.233 |
| | $\infty$ | 0 | 4.323 | 1.794 | 1.374 | 0.955 | 0.955 | 0.955 | 0.955 |
| | | 1 | 2.250 | 1.125 | 0.938 | 0.750 | 0.750 | 0.750 | 0.750 |
| SP | 1 | 0 | 2.333 | 1.126 | 0.935 | 0.760 | 0.705 | 0.724 | 0.774 |
| | | 1 | —— | —— | —— | —— | —— | —— | —— |
| | 2 | 0 | 2.737 | 1.289 | 1.063 | 0.868 | 0.868 | 0.924 | 1.038 |
| | | 1 | 2.301 | 1.230 | 1.067 | 0.934 | 0.965 | 1.042 | 1.168 |
| | 3 | 0 | 2.913 | 1.348 | 1.105 | 0.904 | 0.924 | 0.993 | 1.136 |
| | | 1 | 2.277 | 1.225 | 1.070 | 0.957 | 1.033 | 1.141 | 1.317 |
| | 4 | 0 | 3.016 | 1.378 | 1.125 | 0.920 | 0.949 | 1.020 | 1.170 |
| | | 1 | 2.225 | 1.200 | 1.051 | 0.956 | 1.057 | 1.179 | 1.383 |
| | 6 | 0 | 3.137 | 1.410 | 1.142 | 0.932 | 0.966 | 1.032 | 1.176 |
| | | 1 | 2.128 | 1.146 | 1.007 | 0.936 | 1.057 | 1.189 | 1.414 |
| | 10 | 0 | 3.255 | 1.438 | 1.154 | 0.937 | 0.969 | 1.022 | 1.139 |
| | | 1 | 2.001 | 1.068 | 0.936 | 0.891 | 1.017 | 1.144 | 1.365 |
| | $\infty$ | 0 | 3.507 | 1.503 | 1.176 | 0.895 | 0.895 | 0.895 | 0.895 |
| | | 1 | 1.667 | 0.833 | 0.694 | 0.556 | 0.556 | 0.556 | 0.556 |

**7. Simulations.** We conducted two simulations, one for pure scale alternatives and another one for pure shape alternatives, both in dimension $k=2$. More precisely, starting from two sets of i.i.d. bivariate random vectors $\boldsymbol{\varepsilon}_{1j}$ $(j=1,\ldots,n_1=100)$ and $\boldsymbol{\varepsilon}_{2j}$ $(j=1,\ldots,n_2=100)$ with spherical densities (the standard bivariate normal and bivariate $t$-distributions with 0.5, 2 and 5 degrees of freedom) centered at $\mathbf{0}$, we considered independent samples obtained from

$$\mathbf{X}_{1j}=\mathbf{A}_1\boldsymbol{\varepsilon}_{1j}+\boldsymbol{\theta}_1, \qquad j=1,\ldots,n_1,$$

and

$$\mathbf{X}_{2j}=\mathbf{A}_2(\ell)\boldsymbol{\varepsilon}_{2j}+\boldsymbol{\theta}_2, \qquad j=1,\ldots,n_2,$$

where $\mathbf{A}_2(\ell)\mathbf{A}_2'(\ell)=(1+\ell s^2)(\mathbf{A}_1\mathbf{A}_1'+\ell\mathbf{v})$ [$\mathbf{v}$ a symmetric $(k\times k)$ matrix with $\mathrm{tr}((\mathbf{A}_1\mathbf{A}_1')^{-1}\mathbf{v})=0$], $\ell=0,1,2,3$. The values of $\ell$ produce the null ($\ell=0$) and increasingly heterogeneous alternatives ($\ell=1,2,3$); all tests being affine-invariant, there is no loss of generality in letting $\mathbf{A}_1=\mathbf{I}_2$ and $\boldsymbol{\theta}_1=\boldsymbol{\theta}_2=\mathbf{0}$.

In the first simulation (pure scale alternatives), we generated $N=2{,}500$ independent samples, with $\mathbf{v}=\mathbf{0}$ and $s^2=0.30$, 0.44, 0.56 and 1.50 under Gaussian, $t_5$, $t_2$ and $t_{0.5}$ alternatives, respectively; these values of $s^2$ have been chosen in order to obtain rejection frequencies of the same order under those four densities. In the second simulation (pure shape alternatives), we similarly generated $N=2{,}500$ independent samples, with $s^2=0$ and $\overset{\circ}{\mathrm{vech}}\,\mathbf{v}=(0,0.18)'$, $(0,0.20)'$, $(0,0.21)'$ and $(0,0.22)'$ under Gaussian, $t_5$, $t_2$ and $t_{0.5}$ alternatives, respectively, still with the same objective of obtaining comparable empirical powers under the various densities considered.

In each of these samples, we performed the following eight tests (all at asymptotic level $\alpha=5\%$): (a) the modified likelihood ratio test $\phi^{(n)}_{\mathrm{MLRT}}$; (b) the parametric Gaussian test $\phi^{(n)}_{\mathcal{N}}$ (equivalently, Schott's original test $\phi^{(n)}_{\mathrm{Schott}}$); (c) its pseudo-Gaussian version $\phi^{(n)}_{\mathcal{N}*}$, based on (5.9) (equivalently, the modified Schott test $\phi^{(n)}_{\mathrm{Schott}*}$); (d) the van der Waerden test $\underset{\sim}{\phi}{}^{(n)}_{\mathrm{vdW}}$ [based on (5.6)]; (e)–(g) $t_\nu$-score tests $\underset{\sim}{\phi}{}^{(n)}_{f^t_{1,\nu}}$ with $\nu=5$, 2 and 0.5 [based on (5.7)], as well as (h) the Spearman test $\underset{\sim}{\phi}{}^{(n)}_{K_2}$ [based on (5.8)]. It can be checked that the Wilcoxon test $\underset{\sim}{\phi}{}^{(n)}_{K_1}$, for $k=2$, coincides with $\underset{\sim}{\phi}{}^{(n)}_{f^t_{1,2}}$.

Rejection frequencies are reported in Table 2 for pure scale alternatives, and in Table 3 for pure shape alternatives [the corresponding individual confidence intervals (for $N=2{,}500$ replications), at confidence level 0.95, have

half-widths 0.0044, 0.0080 and 0.0100, for frequencies of the order of 0.05 (0.95), 0.20 (0.80) and 0.50, resp.]. These frequencies indicate that the rank tests, when based on their asymptotic chi-square critical values, are conservative and significantly biased at moderate sample size (100 observations in each group). In order to remedy this, we also implemented versions of each of the rank procedures based on estimations of the (distribution-free) quantile of the test statistic under known location and known common null value of the shape. These estimations, just as the asymptotic chi-square quantile, are consistent approximations of the corresponding exact quantiles under the null, and were obtained for each of the five rank tests under consideration in (d)–(h) above, as the empirical 0.05-upper quantiles $q_{0.95}$ of each rank-based test statistic in a collection of $10^5$ simulated multinormal samples, yielding $q_{0.95} = 7.2117$, 7.6351, 7.7473, 7.7636 and 7.6773, respectively. These bias-corrected critical values all are smaller than the asymptotic chi-square one $\chi^2_{3;0.95} = 7.8147$, so that the resulting tests are uniformly less conservative than the original ones. The resulting rejection frequencies are given in parentheses.

Inspection of Tables 2 and 3 confirms the fact that the parametric Gaussian tests $\phi_{\mathcal{N}}$, contrary to the pseudo-Gaussian ones $\phi_{\mathcal{N}*}$, are invalid under non-Gaussian densities (culminating, under $t_{0.5}$, with a size of 0.9992). However, even the pseudo-Gaussian tests $\phi_{\mathcal{N}*}$, though resisting non-Gaussian densities with finite fourth-order moments, are collapsing under the heavy-tailed $t_{0.5}$ and $t_2$ distributions (with power less than $10^{-4}$ under $t_{0.5}$). In sharp contrast with this, all rank-based tests appear to satisfy the 5% probability level constraint. They are conservative in their original versions (particularly so for van der Waerden scores), but reasonably unbiased (for $n_1 = n_2 = 100$) after bias-correction. Empirical power rankings are essentially consistent with ARE values; in order to allow for meaningful comparisons under infinite fourth-order moments, we also provide AREs with respect to the van der Waerden rank test.

## APPENDIX

### A.1. Proofs of Lemma 4.1, Theorem 5.1 and Corollary 5.1.

PROOF OF LEMMA 4.1. (i) Fix $r \in \{1, \ldots, m\}$. Clearly, under $\mathrm{P}^{(n)}_{\boldsymbol{\vartheta};g_1}$, $\underset{\sim}{\Delta}^{II,r}_{\boldsymbol{\vartheta};K} = \Delta^{II,r}_{\boldsymbol{\vartheta};K;g_1} + o_{L^2}(1)$ iff

$$\text{(A.1)}\quad \sum_{i=1}^{m}\sum_{j=1}^{n_i} c^{(n)}_{ij;r} K\left(\frac{R_{ij}}{n+1}\right) = \sum_{i=1}^{m}\sum_{j=1}^{n_i} c^{(n)}_{ij;r} K\left(\tilde{G}_{1k}\left(\frac{d_{ij}}{\sigma}\right)\right) + o_{L^2}(1),$$

where $c^{(n)}_{ij;r} := n_i^{-1/2}\delta_{i,r}$. For $\boldsymbol{\vartheta} \in \mathcal{M}(\boldsymbol{\Upsilon})$, the $R_{ij}$'s are the ranks of the $d_{ij}/\sigma$'s, which under $\mathrm{P}^{(n)}_{\boldsymbol{\vartheta};g_1}$ are i.i.d. with distribution function $\tilde{G}_{1k}$. The asymptotic

TABLE 2

*Rejection frequencies (out of $N = 2{,}500$ replications), under the null and various scale alternatives (see Section 7 for details), of the Gaussian modified LRT ($\phi_{\mathrm{MLRT}}$), the parametric Gaussian test ($\phi_{\mathcal{N}}$), its pseudo-Gaussian version ($\phi_{\mathcal{N}*}$), and the signed-rank van der Waerden ($\underset{\sim}{\phi}_{\mathrm{vdW}}$), $t_\nu$-score ($\underset{\sim}{\phi}_{f^t_{1,\nu}}$, $\nu = 0.5$, 2, 5), Wilcoxon-type ($\underset{\sim}{\phi}_{K_1}$) and Spearman-type ($\underset{\sim}{\phi}_{K_2}$) tests, respectively. Sample sizes are $n_1 = n_2 = 100$. ARE values are provided with respect to the parametric pseudo-Gaussian ($\mathrm{ARE}_{\mathcal{N}*}$) and van der Waerden rank tests ($\mathrm{ARE}_{\mathrm{vdW}}$); "ND" means "not defined" (this occurs as soon as one the tests involved is not valid under the distribution under consideration)*

| | | $\ell$ | | | | | |
|---|---|---|---|---|---|---|---|
| **Test** | | **0** | **1** | **2** | **3** | $\mathbf{ARE}_{\mathcal{N}*}$ | $\mathbf{ARE}_{\mathbf{vdW}}$ |
| $\phi_{\mathrm{LRT}}$ | $\mathcal{N}$ | 0.0512 | 0.3168 | 0.7932 | 0.9776 | 1.000 | 1.000 |
| $\phi_{\mathrm{MLRT}}$ | | 0.0500 | 0.3100 | 0.7876 | 0.9772 | 1.000 | 1.000 |
| $\phi_{\mathcal{N}}$ | | 0.0464 | 0.3008 | 0.7760 | 0.9756 | 1.000 | 1.000 |
| $\phi_{\mathcal{N}*}$ | | 0.0472 | 0.2944 | 0.7568 | 0.9736 | 1.000 | 1.000 |
| $\underset{\sim}{\phi}_{\mathrm{vdW}}$ | | 0.0348 (0.0472) | 0.2388 (0.2932) | 0.6912 (0.7316) | 0.9520 (0.9676) | 1.000 | 1.000 |
| $\underset{\sim}{\phi}_{f^t_{1,5}}$ | | 0.0444 (0.0496) | 0.2604 (0.2724) | 0.7080 (0.7200) | 0.9552 (0.9600) | 0.918 | 0.918 |
| $\underset{\sim}{\phi}_{f^t_{1,2}} = \underset{\sim}{\phi}_{K_1}$ | | 0.0516 (0.0524) | 0.2180 (0.2248) | 0.6360 (0.6404) | 0.9004 (0.9016) | 0.750 | 0.750 |
| $\underset{\sim}{\phi}_{f^t_{1,0.5}}$ | | 0.0476 (0.0492) | 0.1224 (0.1248) | 0.3252 (0.3260) | 0.5692 (0.5724) | 0.360 | 0.360 |
| $\underset{\sim}{\phi}_{K_2}$ | | 0.0432 (0.0480) | 0.2448 (0.2572) | 0.6956 (0.7060) | 0.9480 (0.9508) | 0.868 | 0.868 |
| $\phi_{\mathrm{LRT}}$ | $t_5$ | 0.3288 | 0.6308 | 0.9168 | 0.9872 | ND | ND |
| $\phi_{\mathrm{MLRT}}$ | | 0.3244 | 0.6260 | 0.9144 | 0.9868 | ND | ND |
| $\phi_{\mathcal{N}}$ | | 0.3160 | 0.6208 | 0.9092 | 0.9856 | ND | ND |
| $\phi_{\mathcal{N}*}$ | | 0.0300 | 0.1896 | 0.5268 | 0.7892 | 1.000 | 0.392 |
| $\underset{\sim}{\phi}_{\mathrm{vdW}}$ | | 0.0320 (0.0444) | 0.2500 (0.2956) | 0.7068 (0.7468) | 0.9396 (0.9560) | 2.551 | 1.000 |
| $\underset{\sim}{\phi}_{f^t_{1,5}}$ | | 0.0428 (0.0480) | 0.3004 (0.3152) | 0.7740 (0.7812) | 0.9636 (0.9676) | 2.778 | 1.089 |
| $\underset{\sim}{\phi}_{f^t_{1,2}} = \underset{\sim}{\phi}_{K_1}$ | | 0.0488 (0.0512) | 0.2916 (0.2980) | 0.7456 (0.7520) | 0.9528 (0.9544) | 2.604 | 1.021 |

Table 2
*Continued*

| **Test** | | $\ell$ **0** | **1** | **2** | **3** | $\mathbf{ARE}_{\mathcal{N}*}$ | $\mathbf{ARE}_{\mathbf{vdW}}$ |
|---|---|---|---|---|---|---|---|
| $\underset{\sim}{\phi}_{f^t_{1,0.5}}$ | | 0.0512 (0.0516) | 0.1824 (0.1848) | 0.4916 (0.4972) | 0.7556 (0.7612) | 1.543 | 0.605 |
| $\underset{\sim}{\phi}_{K_2}$ | | 0.0448 (0.0484) | 0.3068 (0.3184) | 0.7720 (0.7828) | 0.9644 (0.9656) | 2.737 | 1.073 |
| $\phi_{\mathrm{LRT}}$ | $t_2$ | 0.8728 | 0.9164 | 0.9496 | 0.9712 | ND | ND |
| $\phi_{\mathrm{MLRT}}$ | | 0.8696 | 0.9156 | 0.9496 | 0.9700 | ND | ND |
| $\phi_{\mathcal{N}}$ | | 0.8648 | 0.9120 | 0.9480 | 0.9684 | ND | ND |
| $\phi_{\mathcal{N}*}$ | | 0.0120 | 0.0300 | 0.0672 | 0.1276 | ND | ND |
| $\underset{\sim}{\phi}_{\mathrm{vdW}}$ | | 0.0428 (0.0568) | 0.1880 (0.2264) | 0.5368 (0.5816) | 0.7988 (0.8324) | ND | 1.000 |
| $\underset{\sim}{\phi}_{f^t_{1,5}}$ | | 0.0536 (0.0592) | 0.2532 (0.2644) | 0.6592 (0.6704) | 0.9000 (0.9072) | ND | 1.250 |
| $\underset{\sim}{\phi}_{f^t_{1,2}} = \underset{\sim}{\phi}_{K_1}$ | | 0.0508 (0.0532) | 0.2732 (0.2816) | 0.6912 (0.6964) | 0.9212 (0.9236) | ND | 1.333 |
| $\underset{\sim}{\phi}_{f^t_{1,0.5}}$ | | 0.0496 (0.0500) | 0.2116 (0.2136) | 0.5404 (0.5468) | 0.8128 (0.8144) | ND | 1.000 |
| $\underset{\sim}{\phi}_{K_2}$ | | 0.0572 (0.0588) | 0.2568 (0.2652) | 0.6632 (0.6708) | 0.9036 (0.9080) | ND | 1.250 |
| $\phi_{\mathrm{LRT}}$ | $t_{0.5}$ | 0.9992 | 0.9996 | 0.9996 | 0.9988 | ND | ND |
| $\phi_{\mathrm{MLRT}}$ | | 0.9992 | 0.9996 | 0.9996 | 0.9988 | ND | ND |
| $\phi_{\mathcal{N}}$ | | 0.9992 | 0.9996 | 0.9988 | 0.9988 | ND | ND |
| $\phi_{\mathcal{N}*}$ | | 0 | 0 | 0 | 0 | ND | ND |
| $\underset{\sim}{\phi}_{\mathrm{vdW}}$ | | 0.0388 (0.0520) | 0.1464 (0.1764) | 0.3096 (0.3572) | 0.4608 (0.5188) | ND | 1.000 |
| $\underset{\sim}{\phi}_{f^t_{1,5}}$ | | 0.0496 (0.0524) | 0.2328 (0.2452) | 0.5000 (0.5132) | 0.6920 (0.7044) | ND | 1.543 |
| $\underset{\sim}{\phi}_{f^t_{1,2}} = \underset{\sim}{\phi}_{K_1}$ | | 0.0508 (0.0528) | 0.3076 (0.3136) | 0.6404 (0.6448) | 0.8276 (0.8316) | ND | 2.083 |
| $\underset{\sim}{\phi}_{f^t_{1,0.5}}$ | | 0.0604 (0.0616) | 0.3928 (0.3972) | 0.7572 (0.7600) | 0.9208 (0.9212) | ND | 2.778 |
| $\underset{\sim}{\phi}_{K_2}$ | | 0.0488 (0.0524) | 0.2136 (0.2228) | 0.4728 (0.4840) | 0.6672 (0.6792) | ND | 1.435 |

TABLE 3

*Rejection frequencies (out of $N = 2{,}500$ replications), under the null and various shape alternatives (see Section 7 for details), of the Gaussian modified LRT ($\phi_{\mathrm{MLRT}}$), the parametric Gaussian test ($\phi_{\mathcal{N}}$), its pseudo-Gaussian version ($\phi_{\mathcal{N}*}$), and the signed-rank van der Waerden ($\underset{\sim}{\phi}_{\mathrm{vdW}}$), $t_\nu$-score ($\underset{\sim}{\phi}_{f^t_{1,\nu}}$, $\nu = 0.5$, 2, 5), Wilcoxon-type ($\underset{\sim}{\phi}_{K_1}$) and Spearman-type ($\underset{\sim}{\phi}_{K_2}$) tests, respectively. Sample sizes are $n_1 = n_2 = 100$. ARE values are provided with respect to the parametric pseudo-Gaussian ($\mathrm{ARE}_{\mathcal{N}*}$) and van der Waerden rank tests ($\mathrm{ARE}_{\mathrm{vdW}}$); "ND" means "not defined" (this occurs as soon as one the tests involved is not valid under the distribution under consideration)*

| | | $\boldsymbol{\ell}$ | | | | | |
|---|---|---|---|---|---|---|---|
| **Test** | | **0** | **1** | **2** | **3** | $\mathbf{ARE}_{\mathcal{N}*}$ | $\mathbf{ARE}_{\mathbf{vdW}}$ |
| $\phi_{\mathrm{LRT}}$ | $\mathcal{N}$ | 0.0512 | 0.1564 | 0.6032 | 0.9668 | 1.000 | 1.000 |
| $\phi_{\mathrm{MLRT}}$ | | 0.0500 | 0.1532 | 0.5984 | 0.9656 | 1.000 | 1.000 |
| $\phi_{\mathcal{N}}$ | | 0.0464 | 0.1484 | 0.5900 | 0.9640 | 1.000 | 1.000 |
| $\phi_{\mathcal{N}*}$ | | 0.0472 | 0.1444 | 0.5812 | 0.9648 | 1.000 | 1.000 |
| $\underset{\sim}{\phi}_{\mathrm{vdW}}$ | | 0.0348 (0.0472) | 0.1212 (0.1464) | 0.5248 (0.5828) | 0.9488 (0.9632) | 1.000 | 1.000 |
| $\underset{\sim}{\phi}_{f^t_{1,5}}$ | | 0.0444 (0.0496) | 0.1452 (0.1536) | 0.5456 (0.5596) | 0.9464 (0.9496) | 0.945 | 0.945 |
| $\underset{\sim}{\phi}_{f^t_{1,2}} = \underset{\sim}{\phi}_{K_1}$ | | 0.0516 (0.0524) | 0.1364 (0.1392) | 0.4928 (0.5004) | 0.9272 (0.9276) | 0.844 | 0.844 |
| $\underset{\sim}{\phi}_{f^t_{1,0.5}}$ | | 0.0476 (0.0492) | 0.1120 (0.1140) | 0.3996 (0.4036) | 0.8356 (0.8388) | 0.648 | 0.648 |
| $\underset{\sim}{\phi}_{K_2}$ | | 0.0432 (0.0480) | 0.1440 (0.1508) | 0.5420 (0.5512) | 0.9460 (0.9488) | 0.934 | 0.934 |
| $\phi_{\mathrm{LRT}}$ | $t_5$ | 0.3288 | 0.4632 | 0.7840 | 0.9808 | ND | ND |
| $\phi_{\mathrm{MLRT}}$ | | 0.3244 | 0.4600 | 0.7816 | 0.9800 | ND | ND |
| $\phi_{\mathcal{N}}$ | | 0.3160 | 0.4512 | 0.7728 | 0.9796 | ND | ND |
| $\phi_{\mathcal{N}*}$ | | 0.0300 | 0.1020 | 0.4204 | 0.8552 | 1.000 | 0.454 |
| $\underset{\sim}{\phi}_{\mathrm{vdW}}$ | | 0.0320 (0.0444) | 0.1268 (0.1592) | 0.5320 (0.5816) | 0.9576 (0.9692) | 2.204 | 1.000 |
| $\underset{\sim}{\phi}_{f^t_{1,5}}$ | | 0.0428 (0.0480) | 0.1572 (0.1676) | 0.5928 (0.6036) | 0.9720 (0.9740) | 2.333 | 1.059 |
| $\underset{\sim}{\phi}_{f^t_{1,2}} = \underset{\sim}{\phi}_{K_1}$ | | 0.0488 (0.0512) | 0.1608 (0.1632) | 0.5876 (0.5916) | 0.9684 (0.9692) | 2.258 | 1.024 |

TABLE 3
*Continued*

| Test | | $\ell$ = 0 | 1 | 2 | 3 | $\mathbf{ARE}_{\mathcal{N}*}$ | $\mathbf{ARE}_{\mathrm{vdW}}$ |
|---|---|---|---|---|---|---|---|
| $\underset{\sim}{\phi}_{f^t_{1,0.5}}$ | | 0.0512 (0.0516) | 0.1376 (0.1388) | 0.5088 (0.5132) | 0.9312 (0.9332) | 1.896 | 0.860 |
| $\underset{\sim}{\phi}_{K_2}$ | | 0.0448 (0.0484) | 0.1612 (0.1704) | 0.5860 (0.5976) | 0.9700 (0.9716) | 2.301 | 1.044 |
| $\phi_{\mathrm{LRT}}$ | $t_2$ | 0.8728 | 0.8912 | 0.9376 | 0.9768 | ND | ND |
| $\phi_{\mathrm{MLRT}}$ | | 0.8696 | 0.8892 | 0.9364 | 0.9768 | ND | ND |
| $\phi_{\mathcal{N}}$ | | 0.8648 | 0.8864 | 0.9332 | 0.9764 | ND | ND |
| $\phi_{\mathcal{N}*}$ | | 0.0120 | 0.0224 | 0.0808 | 0.2380 | ND | ND |
| $\underset{\sim}{\phi}_{\mathrm{vdW}}$ | | 0.0428 (0.0568) | 0.1180 (0.1488) | 0.4596 (0.5120) | 0.9216 (0.9416) | ND | 1.000 |
| $\underset{\sim}{\phi}_{f^t_{1,5}}$ | | 0.0536 (0.0592) | 0.1488 (0.1560) | 0.5460 (0.5572) | 0.9576 (0.9616) | ND | 1.147 |
| $\underset{\sim}{\phi}_{f^t_{1,2}} = \underset{\sim}{\phi}_{K_1}$ | | 0.0508 (0.0532) | 0.1584 (0.1612) | 0.5640 (0.5668) | 0.9668 (0.9668) | ND | 1.185 |
| $\underset{\sim}{\phi}_{f^t_{1,0.5}}$ | | 0.0496 (0.0500) | 0.1508 (0.1524) | 0.5212 (0.5256) | 0.9412 (0.9420) | ND | 1.089 |
| $\underset{\sim}{\phi}_{K_2}$ | | 0.0572 (0.0588) | 0.1440 (0.1500) | 0.5288 (0.5420) | 0.9516 (0.9564) | ND | 1.111 |
| $\phi_{\mathrm{LRT}}$ | $t_{0.5}$ | 0.9992 | 0.9988 | 0.9992 | 0.9992 | ND | ND |
| $\phi_{\mathrm{MLRT}}$ | | 0.9992 | 0.9988 | 0.9992 | 0.9992 | ND | ND |
| $\phi_{\mathcal{N}}$ | | 0.9992 | 0.9988 | 0.9992 | 0.9992 | ND | ND |
| $\phi_{\mathcal{N}*}$ | | 0 | 0 | 0.0004 | 0.0008 | ND | ND |
| $\underset{\sim}{\phi}_{\mathrm{vdW}}$ | | 0.0388 (0.0520) | 0.0964 (0.1208) | 0.3328 (0.3792) | 0.7960 (0.8328) | ND | 1.000 |
| $\underset{\sim}{\phi}_{f^t_{1,5}}$ | | 0.0496 (0.0524) | 0.1280 (0.1356) | 0.4288 (0.4408) | 0.8928 (0.9004) | ND | 1.254 |
| $\underset{\sim}{\phi}_{f^t_{1,2}} = \underset{\sim}{\phi}_{K_1}$ | | 0.0508 (0.0528) | 0.1396 (0.1440) | 0.4840 (0.4880) | 0.9360 (0.9380) | ND | 1.418 |
| $\underset{\sim}{\phi}_{f^t_{1,0.5}}$ | | 0.0604 (0.0616) | 0.1644 (0.1648) | 0.5356 (0.5388) | 0.9560 (0.9568) | ND | 1.543 |
| $\underset{\sim}{\phi}_{K_2}$ | | 0.0488 (0.0524) | 0.1208 (0.1272) | 0.3968 (0.4064) | 0.8624 (0.8704) | ND | 1.138 |

equivalence (A.1) thus follows from Hájek's classical projection result for linear rank statistics (see, e.g., [32], Chapter 2), since (a) the $c_{ij;r}^{(n)}$'s are not all equal and (b)

$$\frac{\max_{i,j}(c_{ij;r}^{(n)} - n^{-1}\sum_{i,j} c_{ij;r}^{(n)})^2}{\sum_{i,j}(c_{ij;r}^{(n)} - n^{-1}\sum_{i,j} c_{ij;r}^{(n)})^2} = n^{-1}\max\left(\frac{1-\lambda_r^{(n)}}{\lambda_r^{(n)}}, \frac{\lambda_r^{(n)}}{1-\lambda_r^{(n)}}\right) = o(1) \qquad \text{as } n \to \infty$$

(the Noether condition) holds; see the comments after Assumption (B).

Similarly, for the shape part, $\underset{\sim}{\mathbf{\Delta}}_{\boldsymbol{\vartheta};K}^{III,r} = \mathbf{\Delta}_{\boldsymbol{\vartheta};K;g_1}^{III,r} + o_{L^2}(1)$ under $\mathrm{P}_{\boldsymbol{\vartheta};g_1}^{(n)}$ iff

$$n_r^{-1/2}\mathbf{M}_k(\mathbf{V})(\mathbf{V}^{\otimes 2})^{-1/2}\sum_{j=1}^{n_r}\left[K\left(\frac{R_{rj}}{n+1}\right) - K\left(\tilde{G}_{1k}\left(\frac{d_{rj}}{\sigma}\right)\right)\right]\mathbf{J}_k^{\perp}\operatorname{vec}(\mathbf{U}_{rj}\mathbf{U}_{rj}') = o_{L^2}(1)$$

[where $\mathbf{J}_k^{\perp} := \mathbf{I}_{k^2} - \frac{1}{k}\mathbf{J}_k$ satisfies $\mathbf{M}_k(\mathbf{V})(\mathbf{V}^{\otimes 2})^{-1/2}\mathbf{J}_k^{\perp} = \mathbf{M}_k(\mathbf{V})(\mathbf{V}^{\otimes 2})^{-1/2}$ and is such that $\mathbf{J}_k^{\perp}\operatorname{vec}(\mathbf{U}_{rj}\mathbf{U}_{rj}')$ is exactly centered], or equivalently iff

$$T_{r;l}^{(n)} := \sum_{i=1}^{m}\sum_{j=1}^{n_i} c_{ij;r}^{(n)}\left[K\left(\frac{R_{ij}}{n+1}\right) - K\left(\tilde{G}_{1k}\left(\frac{d_{ij}}{\sigma}\right)\right)\right][\mathbf{J}_k^{\perp}\operatorname{vec}(\mathbf{U}_{ij}\mathbf{U}_{ij}')]_\ell = o_{L^2}(1), \tag{A.2}$$

for all $\ell \in \{1, 2, \ldots, k^2\}$, still under $\mathrm{P}_{\boldsymbol{\vartheta};g_1}^{(n)}$. Now,

$$\mathrm{E}[(T_{r;\ell}^{(n)})^2] = C_{\ell,k}\sum_{i=1}^{m}\sum_{j=1}^{n_i}(c_{ij;r}^{(n)})^2\mathrm{E}\left[\left(K\left(\frac{R_i}{n+1}\right) - K\left(\tilde{G}_{1k}\left(\frac{d_i}{\sigma}\right)\right)\right)^2\right]$$

where, denoting by $U_{ij,s}$ the $s$th component of $\mathbf{U}_{ij}$, $C_{\ell,k} = \mathrm{Var}[U_{11,1}^2] = 2(k-1)/(k^2(k+2))$ for $\ell \in \mathfrak{L}_k := \{mk+m+1, m = 0, 1, \ldots, k-1\}$ and $C_{\ell,k} = \mathrm{Var}[U_{11,1}U_{11,2}] = 1/k^2$ for $\ell \notin \mathfrak{L}_k$. Here, the Hájek projection result for linear *signed*-rank statistics (see, e.g., [32], Chapter 3) yields (A.2), since $\max_{i,j}(c_{ij;r}^{(n)})^2/\sum_{i,j}(c_{ij;r}^{(n)})^2 = n_r^{-1} = o(1)$, as $n \to \infty$.

As for (ii), the result straightforwardly follows, under $\mathrm{P}_{\boldsymbol{\vartheta};g_1}^{(n)}$ with $\boldsymbol{\vartheta} \in \mathcal{M}(\boldsymbol{\Upsilon})$, from the multivariate CLT. The result under local alternatives is obtained as usual, by establishing the joint normality under $\mathrm{P}_{\boldsymbol{\vartheta};g_1}^{(n)}$ of $\mathbf{\Delta}_{\boldsymbol{\vartheta};K;g_1}$ and $\Lambda_{\boldsymbol{\vartheta}+n^{-1/2}\boldsymbol{\nu}^{(n)}\boldsymbol{\tau}/\boldsymbol{\vartheta};g_1}^{(n)}$, then applying Le Cam's third lemma; the required joint normality follows from a routine application of the classical Cramér–Wold device. □

Proof of Theorem 5.1. (i) The continuity of the mapping $\boldsymbol{\vartheta} \mapsto (\mathbf{P}^{II}_{\boldsymbol{\vartheta};K}, \mathbf{P}^{III}_{\boldsymbol{\vartheta};K})$, Proposition 4.2 (jointly with Assumption (D1) and the fact that $[\mathbf{I}_m - \mathbf{C}^{(n)}](\boldsymbol{\Lambda}^{(n)})^{-1}\mathbf{1}_m = \mathbf{0})$, and Lemma 4.1(i), entail

$$\begin{aligned} \underset{\sim}{Q}^{(n)}_K &= (\underset{\sim}{\boldsymbol{\Delta}}^{II}_{\boldsymbol{\vartheta};K})'\mathbf{P}^{II}_{\boldsymbol{\vartheta};K}\underset{\sim}{\boldsymbol{\Delta}}^{II}_{\boldsymbol{\vartheta};K} + (\underset{\sim}{\boldsymbol{\Delta}}^{III}_{\boldsymbol{\vartheta};K})'\mathbf{P}^{III}_{\boldsymbol{\vartheta};K}\underset{\sim}{\boldsymbol{\Delta}}^{III}_{\boldsymbol{\vartheta};K} + o_{\mathrm{P}}(1) \\ &= (\boldsymbol{\Delta}^{II}_{\boldsymbol{\vartheta};K;g_1})'\mathbf{P}^{II}_{\boldsymbol{\vartheta};K}\boldsymbol{\Delta}^{II}_{\boldsymbol{\vartheta};K;g_1} + (\boldsymbol{\Delta}^{III}_{\boldsymbol{\vartheta};K;g_1})'\mathbf{P}^{III}_{\boldsymbol{\vartheta};K}\boldsymbol{\Delta}^{III}_{\boldsymbol{\vartheta};K;g_1} + o_{\mathrm{P}}(1) \end{aligned} \tag{A.3}$$

under $\mathrm{P}^{(n)}_{\boldsymbol{\vartheta};g_1}$, $\boldsymbol{\vartheta} \in \mathcal{M}(\boldsymbol{\Upsilon})$ (and therefore, also under the contiguous sequence $\mathrm{P}^{(n)}_{\boldsymbol{\vartheta}+n^{-1/2}\boldsymbol{\nu}^{(n)}\boldsymbol{\tau}^{(n)};g_1}$). Now, since $(\boldsymbol{\Gamma}^{II}_{\boldsymbol{\vartheta};K})^{1/2}\mathbf{P}^{II}_{\boldsymbol{\vartheta};K}(\boldsymbol{\Gamma}^{II}_{\boldsymbol{\vartheta};K})^{1/2}$ is a symmetric idempotent matrix with rank $m-1$, it follows from Lemma 4.1(ii) that the first term in (A.3) is asymptotically chi-square with $m-1$ degrees of freedom under $\mathrm{P}^{(n)}_{\boldsymbol{\vartheta};g_1}$, $\boldsymbol{\vartheta} \in \mathcal{M}(\boldsymbol{\Upsilon})$, and asymptotically noncentral chi-square, still with $m-1$ degrees of freedom, but with noncentrality parameter

$$\left(\frac{\mathcal{L}_k(K, g_1)}{4\sigma^4}\right)^2 \lim_{n\to\infty}\{(\boldsymbol{\tau}^{(n)}_{II})'\mathbf{P}^{II}_{\boldsymbol{\vartheta};K}\boldsymbol{\tau}^{(n)}_{II}\} \tag{A.4}$$

under $\mathrm{P}^{(n)}_{\boldsymbol{\vartheta}+n^{-1/2}\boldsymbol{\nu}^{(n)}\boldsymbol{\tau}^{(n)};g_1}$. Evaluation of the limit in (A.4) yields the first term in (5.4).

As for the shape part, using again Lemma 4.1(ii) and the fact that $(\boldsymbol{\Gamma}^{III}_{\boldsymbol{\vartheta};K})^{1/2} \times \mathbf{P}^{III}_{\boldsymbol{\vartheta};K}(\boldsymbol{\Gamma}^{III}_{\boldsymbol{\vartheta};K})^{1/2}$ is symmetric and idempotent with rank $k_0(m-1)$, we obtain similarly that the second term in (A.3) is asymptotically chi-square with $k_0(m-1)$ degrees of freedom under $\mathrm{P}^{(n)}_{\boldsymbol{\vartheta};g_1}$, $\boldsymbol{\vartheta} \in \mathcal{M}(\boldsymbol{\Upsilon})$, and asymptotically noncentral chi-square, still with $k_0(m-1)$ degrees of freedom but with noncentrality parameter

$$(\mathcal{J}_k(K, g_1))^2 \lim_{n\to\infty}\{(\boldsymbol{\tau}^{(n)}_{III})'[\mathbf{I}_m \otimes \mathbf{H}_k(\mathbf{V})]\mathbf{P}^{III}_{\boldsymbol{\vartheta};K}[\mathbf{I}_m \otimes \mathbf{H}_k(\mathbf{V})]\boldsymbol{\tau}^{(n)}_{III}\} \tag{A.5}$$

under $\mathrm{P}^{(n)}_{\boldsymbol{\vartheta}+n^{-1/2}\boldsymbol{\nu}^{(n)}\boldsymbol{\tau}^{(n)};g_1}$. Evaluation of the limit in (A.5) yields the second term in (5.4). As the two terms in (A.3) are asymptotically uncorrelated [see Lemma 4.1(ii) again], they can indeed be treated separately.

(ii) The fact that $\underset{\sim}{\phi}^{(n)}_K$ has asymptotic level $\alpha$ directly follows from the asymptotic null distribution in part (i) and the classical Helly–Bray theorem.

(iii) Optimality is a consequence of the asymptotic equivalence (A.3), under $g_1 = f_1 (\in \mathcal{F}_a)$, of $\underset{\sim}{Q}^{(n)}_{f_1}$ and the locally asymptotically optimal test statistic $Q_{\boldsymbol{\Upsilon}}$, as described in (4.3). □

Proof of Corollary 5.1. (i) Fix $g_1 \in \mathcal{F}_a$, with $\mathcal{L}_k(K, g_1) \neq 0 \neq \mathcal{J}_k(K, g_1)$. Clearly, $\underset{\sim}{\phi}^{(n)}_K$ is consistent under $\mathrm{P}^{(n)}_{\boldsymbol{\vartheta}+n^{-1/2}\boldsymbol{\nu}^{(n)}\boldsymbol{\tau}^{(n)};g_1}$, $\boldsymbol{\vartheta} \in \mathcal{M}(\boldsymbol{\Upsilon})$ iff the corresponding noncentrality parameter in (5.4) is nonzero. Assume the latter is zero. Then, the assumptions on $g_1$ imply that $s_i^2/\sqrt{\lambda_i} = s_{i'}^2/\sqrt{\lambda_{i'}}$ and

$$\mathrm{tr}\left[\left(\mathbf{V}^{-1/2}\left(\frac{\mathbf{v}_i}{\sqrt{\lambda_i}} - \frac{\mathbf{v}_{i'}}{\sqrt{\lambda_{i'}}}\right)\mathbf{V}^{-1/2}\right)^2\right], \tag{A.6}$$

for all $(i, i')$. Now, since $\mathrm{tr}(\mathbf{A}^2) = 0$ implies that $\mathbf{A} = \mathbf{0}$ for any symmetric $k \times k$ matrix $\mathbf{A}$, it follows from (A.6) that $\mathbf{v}_i/\sqrt{\lambda_i} = \mathbf{v}_{i'}/\sqrt{\lambda_{i'}}$ for all $(i, i')$. This is possible only for $\boldsymbol{\nu\tau} \in \mathcal{M}(\boldsymbol{\Upsilon})$, which establishes the result.

(ii) Going back to the definition of $g_1 \mapsto \mathcal{J}_k(K, g_1)$, we have

$$\begin{aligned}\mathcal{J}_k(K, g_1) &= \int_0^\infty K(\tilde{G}_{1k}(r)) r\varphi_{g_1}(r)\tilde{g}_{1k}(r)\,dr \\ &= \frac{1}{\mu_{k-1;g_1}}\int_0^\infty K(\tilde{G}_{1k}(r))(-\dot{g}_1(r))r^k\,dr.\end{aligned}$$

Integrating by parts, $\mathcal{J}_k(K, g_1) = \int_0^\infty [kK(\tilde{G}_{1k}(r)) + \dot{K}'(\tilde{G}_{1k}(r))r\tilde{g}_{1k}(r)] \times \tilde{g}_{1k}(r)\,dr = k^2 + \int_0^\infty \dot{K}'(\tilde{G}_{1k}(r))r(\tilde{g}_{1k}(r))^2\,dr$, so that $\int_0^\infty \dot{K}'(\tilde{G}_{1k}(r)) \times r(\tilde{g}_{1k}(r))^2\,dr > 0$ guarantees that $\mathcal{L}_k(K, g_1) = \mathcal{J}_k(K, g_1) - k^2 > 0$. The result follows from part (i) of the corollary. □

**A.2. Proof of Proposition 4.2.** Consider an arbitrary value $\boldsymbol{\vartheta} = (\boldsymbol{\vartheta}_I', \boldsymbol{\vartheta}_{II}', \boldsymbol{\vartheta}_{III}')' = (\boldsymbol{\theta}_1', \dots, \boldsymbol{\theta}_m', \sigma^2\mathbf{1}_m', \mathbf{1}_m' \otimes (\mathrm{v\overset{\circ}{e}ch}\,\mathbf{V})')' \in \mathcal{M}(\boldsymbol{\Upsilon})$ of the parameter and a (bounded) sequence of corresponding local perturbations $\boldsymbol{\vartheta}^{(n)} := \boldsymbol{\vartheta} + n^{-1/2}\boldsymbol{\nu}^{(n)} \times \boldsymbol{\tau}^{(n)}$, where

$$\begin{aligned}\boldsymbol{\tau}^{(n)} &= (\boldsymbol{\tau}_I^{(n)\prime}, \boldsymbol{\tau}_{II}^{(n)\prime}, \boldsymbol{\tau}_{III}^{(n)\prime})' \\ &= (\mathbf{t}_1^{(n)\prime}, \dots, \mathbf{t}_m^{(n)\prime}, s_1^{2(n)}, \dots, s_m^{2(n)}, (\mathrm{v\overset{\circ}{e}ch}\,\mathbf{v}_1^{(n)})', \dots, (\mathrm{v\overset{\circ}{e}ch}\,\mathbf{v}_m^{(n)})')'\end{aligned}$$

is such that $\boldsymbol{\vartheta}^{(n)} \in \mathcal{M}(\boldsymbol{\Upsilon})$ for all $n$. To prove Proposition 4.2, it is sufficient to show that, under $\mathrm{P}_{\boldsymbol{\vartheta};g_1}^{(n)}$ (where $g_1$ is as in Proposition 4.2),

$$\begin{aligned}&\underset{\sim}{\boldsymbol{\Delta}}{}_{\boldsymbol{\vartheta}^{(n)};K}^{II} - \underset{\sim}{\boldsymbol{\Delta}}{}_{\boldsymbol{\vartheta};K}^{II} + \frac{\mathcal{L}_k(K, g_1)}{4\sigma^4}\boldsymbol{\tau}_{II}^{(n)} \quad \text{and} \\ &\underset{\sim}{\boldsymbol{\Delta}}{}_{\boldsymbol{\vartheta}^{(n)};K}^{III} - \underset{\sim}{\boldsymbol{\Delta}}{}_{\boldsymbol{\vartheta};K}^{III} + \mathcal{J}_k(K, g_1)[\mathbf{I}_m \otimes \mathbf{H}_k(\mathbf{V})]\boldsymbol{\tau}_{III}^{(n)}\end{aligned} \tag{A.7}$$

are $o_{\mathrm{P}}(1)$ as $n \to \infty$, since the local asymptotic discreteness of $\hat{\boldsymbol{\vartheta}}$ (see, e.g., [25], Lemma 4.4) allows to replace the nonrandom quantity $\boldsymbol{\vartheta}^{(n)}$ with the random one $\hat{\boldsymbol{\vartheta}}$ in (A.7). Note that $\hat{\boldsymbol{\vartheta}}$ being constrained allows us to restrict to $\boldsymbol{\vartheta}^{(n)} \in \mathcal{M}(\boldsymbol{\Upsilon})$. Looking at block $i$ ($i \in \{1, \dots, m\}$), Proposition 4.2 thus is a corollary of the following result.

PROPOSITION A.1. *Assume that* (A), (B) *and* (C) *hold, and that* $g_1 \in \mathcal{F}_a$. *Fix* $\boldsymbol{\vartheta} \in \mathcal{M}(\boldsymbol{\Upsilon})$ *and a sequence* $\boldsymbol{\vartheta}^{(n)} \in \mathcal{M}(\boldsymbol{\Upsilon})$ *as above. Then, for all* $i = 1, \dots, m$,

$$\underset{\sim}{\Delta}^{II,i}_{\boldsymbol{\vartheta}^{(n)};K} - \underset{\sim}{\Delta}^{II,i}_{\boldsymbol{\vartheta};K} + \frac{\mathcal{L}_k(K,g_1)}{4\sigma^4} s_i^{2(n)} \quad \textit{and} \tag{A.8}$$

$$\underset{\sim}{\boldsymbol{\Delta}}^{III,i}_{\boldsymbol{\vartheta}^{(n)};K} - \underset{\sim}{\boldsymbol{\Delta}}^{III,i}_{\boldsymbol{\vartheta};K} + \mathcal{J}_k(K,g_1)\mathbf{H}_k(\mathbf{V})(\overset{\circ}{\mathrm{vech}}\,\mathbf{v}_i^{(n)})$$

*are* $o_{\mathrm{P}}(1)$ *under* $\mathrm{P}^{(n)}_{\boldsymbol{\vartheta};g_1}$, *as* $n \to \infty$.

PROOF. In this proof, we let $\boldsymbol{\theta}_i^n := \boldsymbol{\theta}_i + n_i^{-1/2}\mathbf{t}_i^{(n)}$, $\mathbf{V}^n := \mathbf{V} + n_i^{-1/2}\mathbf{v}_i^{(n)}$, and $\sigma_n^2 := \sigma^2 + n_i^{-1/2}s_i^{2(n)}$ [since $\boldsymbol{\vartheta}, \boldsymbol{\vartheta}^{(n)} \in \mathcal{M}(\boldsymbol{\Upsilon})$, $\sigma_n^2$ and $\mathbf{V}^n$ do not depend on $i$, which explains the notation]. Accordingly, let $\mathbf{Z}_{ij}^0 := \mathbf{V}^{-1/2}(\mathbf{X}_{ij} - \boldsymbol{\theta}_i)$, $d_{ij}^0 := \|\mathbf{Z}_{ij}^0\|$, $\mathbf{U}_{ij}^0 := \mathbf{Z}_{ij}^0/d_{ij}^0$, $\mathbf{Z}_{ij}^n := (\mathbf{V}^n)^{-1/2}(\mathbf{X}_{ij} - \boldsymbol{\theta}_i^n)$, $d_{ij}^n := \|\mathbf{Z}_{ij}^n\|$, and $\mathbf{U}_{ij}^n := \mathbf{Z}_{ij}^n/d_{ij}^n$. Following an argument that goes back to [24], consider the following truncation of the score function $K$: for all $\ell \in \mathbb{N}_0$, define

$$\begin{aligned} K^{(\ell)}(u) := {} & K\left(\frac{2}{\ell}\right)\ell\left(u - \frac{1}{\ell}\right)I_{[1/\ell<u\le 2/\ell]} + K(u)I_{[2/\ell<u\le 1-2/\ell]} \\ & + K\left(1-\frac{2}{\ell}\right)\ell\left(\left(1-\frac{1}{\ell}\right) - u\right)I_{[1-2/\ell<u\le 1-1/\ell]}, \end{aligned}$$

where $I_A$ denotes the indicator of $A$. Since $u \mapsto K(u)$ is continuous, the functions $u \mapsto K^{(\ell)}(u)$ are also continuous on $(0,1)$. It follows that the truncated scores $K^{(\ell)}$ are bounded for all $\ell$. Clearly, it can safely be assumed that $K$ is monotone increasing (rather than the difference of two monotone increasing functions), so that there exists some $L$ such that $|K^{(\ell)}(u)| \le |K(u)|$ for all $u \in (0,1)$ and all $\ell \ge L$.

We start with the proof that the scale part of (A.8) is $o_{\mathrm{P}}(1)$ under $\mathrm{P}^{(n)}_{\boldsymbol{\vartheta};g_1}$. For the shape part, the result is an easy extension of the corresponding result in [14]; details are left to the reader. Lemma 4.1(i) shows that $\underset{\sim}{\Delta}^{II,i}_{\boldsymbol{\vartheta};K} - \Delta^{II,i}_{\boldsymbol{\vartheta};K;g_1}$ is $o_{\mathrm{P}}(1)$, under $\mathrm{P}^{(n)}_{\boldsymbol{\vartheta};g_1}$. Similarly, $\underset{\sim}{\Delta}^{II,i}_{\boldsymbol{\vartheta}^{(n)};K} - \Delta^{II,i}_{\boldsymbol{\vartheta}^{(n)};K;g_1}$ is $o_{\mathrm{P}}(1)$ under $\mathrm{P}^{(n)}_{\boldsymbol{\vartheta}^{(n)};g_1}$—hence, from contiguity, also under $\mathrm{P}^{(n)}_{\boldsymbol{\vartheta};g_1}$. Consequently, (A.8) is asymptotically equivalent, under $\mathrm{P}^{(n)}_{\boldsymbol{\vartheta};g_1}$, to

$$\Delta^{II,i}_{\boldsymbol{\vartheta}^{(n)};K;g_1} - \Delta^{II,i}_{\boldsymbol{\vartheta};K;g_1} + \frac{\mathcal{L}_k(K,g_1)}{4\sigma^4}s_i^{2(n)}. \tag{A.9}$$

Now, $\frac{1}{2}n_i^{-1/2}\sum_{j=1}^{n_i}(K(\tilde{G}_{1k}(d_{ij}^n/\sigma_n)) - k)$, under $\mathrm{P}^{(n)}_{\boldsymbol{\vartheta}^{(n)};g_1}$, is asymptotically normal as $n \to \infty$, with mean zero and variance $\frac{1}{4}\mathcal{L}_k(K)$, so that $\frac{1}{2}(\frac{1}{\sigma_n^2} -$

$\frac{1}{\sigma^2})n_i^{-1/2}\sum_{j=1}^{n_i}(K(\tilde{G}_{1k}(d_{ij}^n/\sigma_n))-k)$ is $o_{\mathrm{P}}(1)$, as $n\to\infty$, under $\mathrm{P}^{(n)}_{\boldsymbol{\vartheta}^{(n)};g_1}$, as well as under $\mathrm{P}^{(n)}_{\boldsymbol{\vartheta};g_1}$ (from contiguity). Consequently, (A.9) is asymptotically equivalent, under $\mathrm{P}^{(n)}_{\boldsymbol{\vartheta};g_1}$, to

$$\mathbf{C}_i^{(n)} := \frac{1}{2\sigma^2}n_i^{-1/2}\sum_{j=1}^{n_i}(K(\tilde{G}_{1k}(d_{ij}^n/\sigma_n))-k) - \frac{1}{2\sigma^2}n_i^{-1/2}\sum_{j=1}^{n_i}(K(\tilde{G}_{1k}(d_{ij}^0/\sigma))-k) + \frac{\mathcal{L}_k(K,g_1)}{4\sigma^4}s_i^{2(n)}.$$

Decompose $\mathbf{C}_i^{(n)}$ into $\mathbf{C}_i^{(n)} = \mathbf{D}_{i1}^{(n;\ell)} + \mathbf{D}_{i2}^{(n;\ell)} - \mathbf{R}_{i1}^{(n;\ell)} + \mathbf{R}_{i2}^{(n;\ell)} + \mathbf{R}_{i3}^{(n;\ell)}$ where, denoting by $\mathrm{E}_0$ expectation under $\mathrm{P}^{(n)}_{\boldsymbol{\vartheta};g_1}$,

$$\begin{aligned}\mathbf{D}_{i1}^{(n;\ell)} &:= \frac{1}{2\sigma^2}n_i^{-1/2}\sum_{j=1}^{n_i}(K^{(\ell)}(\tilde{G}_{1k}(d_{ij}^n/\sigma_n)) - \mathrm{E}[K^{(\ell)}(U)])\\ &\quad - \frac{1}{2\sigma^2}n_i^{-1/2}\sum_{j=1}^{n_i}(K^{(\ell)}(\tilde{G}_{1k}(d_{ij}^0/\sigma)) - \mathrm{E}[K^{(\ell)}(U)])\\ &\quad - \frac{1}{2\sigma^2}n_i^{-1/2}\mathrm{E}_0\left[\sum_{j=1}^{n_i}(K^{(\ell)}(\tilde{G}_{1k}(d_{ij}^n/\sigma_n)) - \mathrm{E}[K^{(\ell)}(U)])\right],\\ \mathbf{D}_{i2}^{(n;\ell)} &:= \frac{1}{2\sigma^2}n_i^{-1/2}\mathrm{E}_0\left[\sum_{j=1}^{n_i}(K^{(\ell)}(\tilde{G}_{1k}(d_{ij}^n/\sigma_n)) - \mathrm{E}[K^{(\ell)}(U)])\right]\\ &\quad + \frac{\mathcal{L}_k(K^{(\ell)},g_1)}{4\sigma^4}s_i^{2(n)},\\ \mathbf{R}_{i1}^{(n;\ell)} &:= \frac{1}{2\sigma^2}n_i^{-1/2}\sum_{j=1}^{n_i}[(K(\tilde{G}_{1k}(d_{ij}^0/\sigma))-k)\\ &\qquad - (K^{(\ell)}(\tilde{G}_{1k}(d_{ij}^0/\sigma)) - \mathrm{E}[K^{(\ell)}(U)])],\\ \mathbf{R}_{i2}^{(n;\ell)} &:= \frac{1}{2\sigma^2}n_i^{-1/2}\sum_{j=1}^{n_i}[(K(\tilde{G}_{1k}(d_{ij}^n/\sigma_n))-k)\\ &\qquad - (K^{(\ell)}(\tilde{G}_{1k}(d_{ij}^n/\sigma_n)) - \mathrm{E}[K^{(\ell)}(U)])],\end{aligned}$$

and

$$\mathbf{R}_{i3}^{(n;\ell)} := \frac{1}{4\sigma^4}(\mathcal{L}_k(K,g_1) - \mathcal{L}_k(K^{(\ell)},g_1))s_i^{2(n)}.$$

We prove that $\mathbf{C}_i^{(n)} = o_{\mathrm{P}}(1)$ [thus completing the proof that (A.8) is $o_{\mathrm{P}}(1)$ under $\mathrm{P}^{(n)}_{\boldsymbol{\vartheta};g_1}$] by establishing that $\mathbf{D}_{i1}^{(n;\ell)}$ and $\mathbf{D}_{i2}^{(n;\ell)}$ are $o_{\mathrm{P}}(1)$ under $\mathrm{P}^{(n)}_{\boldsymbol{\vartheta};g_1}$,

as $n\to\infty$, for fixed $\ell$, and that $\mathbf{R}_{i1}^{(n;\ell)}$, $\mathbf{R}_{i2}^{(n;\ell)}$ and $\mathbf{R}_{i3}^{(n;\ell)}$ are $o_{\mathrm{P}}(1)$ under the same sequence of hypotheses, as $\ell\to\infty$, uniformly in $n$. For the sake of convenience, these three results are treated separately (Lemmas A.1, A.2 and A.3).

LEMMA A.1. *For any fixed $\ell$, $\mathrm{E}_0[|\mathbf{D}_{i1}^{(n;\ell)}|^2]=o(1)$ as $n\to\infty$.*

LEMMA A.2. *For any fixed $\ell$, $\mathbf{D}_{i2}^{(n;\ell)}=o(1)$ as $n\to\infty$.*

LEMMA A.3. *As $\ell\to\infty$, uniformly in $n$,* (i) *$\mathbf{R}_{i1}^{(n;\ell)}$ is $o_{\mathrm{P}}(1)$ under $\mathrm{P}_{\boldsymbol{\vartheta};g_1}^{(n)}$,* (ii) *$\mathbf{R}_{i2}^{(n;\ell)}$ is $o_{\mathrm{P}}(1)$ under $\mathrm{P}_{\boldsymbol{\vartheta};g_1}^{(n)}$ for $n$ sufficiently large and* (iii) *$\mathbf{R}_{i3}^{(n;\ell)}$ is $o(1)$.*

PROOF OF LEMMA A.1. First note that $\mathbf{D}_{i1}^{(n;\ell)}=\frac{1}{2\sigma^2}n_i^{-1/2}\sum_{j=1}^{n_i}[\mathbf{T}_{ij}^{(n;\ell)}-\mathrm{E}_0[\mathbf{T}_{ij}^{(n;\ell)}]]$, where $\mathbf{T}_{ij}^{(n;\ell)}:=K^{(\ell)}(\tilde{G}_{1k}(d_{ij}^{n}/\sigma_n))-K^{(\ell)}(\tilde{G}_{1k}(d_{ij}^{0}/\sigma))$, $j=1,\ldots,n_i$ are i.i.d. Writing $\mathrm{Var}_0$ for variances under $\mathrm{P}_{\boldsymbol{\vartheta};g_1}^{(n)}$,

$$\mathrm{E}_0[|\mathbf{D}_{i1}^{(n;\ell)}|^2]=\mathrm{Var}_0[\mathbf{D}_{i1}^{(n;\ell)}]=\frac{1}{4\sigma^4}\,\mathrm{Var}_0[\mathbf{T}_{i1}^{(n;\ell)}]\le\frac{1}{4\sigma^4}\mathrm{E}_0[|\mathbf{T}_{i1}^{(n;\ell)}|^2],$$

and it only remains to show that, as $n\to\infty$,

$$\begin{aligned}\mathrm{E}_0[|\mathbf{T}_{i1}^{(n;\ell)}|^2]&=\mathrm{E}_0[(K^{(\ell)}(\tilde{G}_{1k}(d_{i1}^{n}/\sigma_n))-K^{(\ell)}(\tilde{G}_{1k}(d_{i1}^{0}/\sigma)))^2]\\&=o(1)\end{aligned}\tag{A.10}$$

Now, $|d_{i1}^{n}/\sigma_n-d_{i1}^{0}/\sigma|\le|d_{i1}^{n}-d_{i1}^{0}|/\sigma_n+|\sigma_n^{-1}-\sigma^{-1}|d_{i1}^{0}$ is $o_{\mathrm{P}}(1)$ under $\mathrm{P}_{\boldsymbol{\vartheta};g_1}^{(n)}$ since $|d_{i1}^{n}-d_{i1}^{0}|$ is; see Lemma A.1 in [14]. This and the continuity of $K^{(\ell)}\circ\tilde{G}_{1k}$ imply that $K^{(\ell)}(\tilde{G}_{1k}(d_{i1}^{n}/\sigma_n))-K^{(\ell)}(\tilde{G}_{1k}(d_{i1}^{0}/\sigma))=o_{\mathrm{P}}(1)$ under $\mathrm{P}_{\boldsymbol{\vartheta};g_1}^{(n)}$, as $n\to\infty$. Since $K^{(\ell)}$ is bounded, this convergence also holds in quadratic mean, which entails (A.10). □

PROOF OF LEMMA A.2. Letting

$$\mathbf{B}_{i1}^{(n;\ell)}:=\frac{1}{2\sigma^2}n_i^{-1/2}\sum_{j=1}^{n_i}(K^{(\ell)}(\tilde{G}_{1k}(d_{ij}^{0}/\sigma))-\mathrm{E}[K^{(\ell)}(U)])$$

one can show that, under $\mathrm{P}_{\boldsymbol{\vartheta};g_1}^{(n)}$, as $n\to\infty$,

$$\mathbf{B}_{i1}^{(n;\ell)}\xrightarrow{\mathcal{L}}\mathcal{N}\Big(\mathbf{0},\frac{1}{4\sigma^4}\,\mathrm{Var}[K^{(\ell)}(U)]\Big).\tag{A.11}$$

Under the sequence of local alternatives $\mathrm{P}^{(n)}_{\boldsymbol{\vartheta}^{(n)};g_1}$, as $n \to \infty$,

$$\mathbf{B}_{i1}^{(n;\ell)} - \frac{\mathcal{L}_k(K^{(\ell)}, g_1)}{4\sigma^4} s_i^{2(n)} \xrightarrow{\mathcal{L}} \mathcal{N}\left(\mathbf{0}, \frac{1}{4\sigma^4}\mathcal{L}_k(K^{(\ell)})\right).$$

Defining $\mathbf{B}_{i2}^{(n;\ell)} := \frac{1}{2\sigma^2} n_i^{-1/2} \sum_{j=1}^{n_i} (K^{(\ell)}(\tilde{G}_{1k}(d_{ij}^n/\sigma_n)) - \mathrm{E}[K^{(\ell)}(U)])$, it follows from ULAN that, under $\mathrm{P}^{(n)}_{\boldsymbol{\vartheta};g_1}$, as $n \to \infty$,

$$\mathbf{B}_{i2}^{(n;\ell)} + \frac{\mathcal{L}_k(K^{(\ell)}, g_1)}{4\sigma^4} s_i^{2(n)} \xrightarrow{\mathcal{L}} \mathcal{N}\left(\mathbf{0}, \frac{1}{4\sigma^4}\mathcal{L}_k(K^{(\ell)})\right). \tag{A.12}$$

Now, from (A.11) and the fact that, under $\mathrm{P}^{(n)}_{\boldsymbol{\vartheta};g_1}$, $\mathbf{D}_{i1}^{(n;\ell)} = \mathbf{B}_{i2}^{(n;\ell)} - \mathbf{B}_{i1}^{(n;\ell)} - \mathrm{E}_0[\mathbf{B}_{i2}^{(n;\ell)}] = o_{\mathrm{P}}(1)$ (Lemma A.1), we obtain that, under $\mathrm{P}^{(n)}_{\boldsymbol{\vartheta};g_1}$, as $n \to \infty$,

$$\mathbf{B}_{i2}^{(n;\ell)} - \mathrm{E}_0[\mathbf{B}_{i2}^{(n;\ell)}] \xrightarrow{\mathcal{L}} \mathcal{N}\left(\mathbf{0}, \frac{1}{4\sigma^4}\mathcal{L}_k(K^{(\ell)})\right). \tag{A.13}$$

The result then follows from comparing (A.12) and (A.13). $\square$

We now complete the proof that (A.8) is $o_{\mathrm{P}}(1)$ under $\mathrm{P}^{(n)}_{\boldsymbol{\vartheta};g_1}$ by proving Lemma A.3.

Proof of Lemma A.3. (i) In view of the independence under $\mathrm{P}^{(n)}_{\boldsymbol{\vartheta};g_1}$ of the $d_{ij}^0$'s,

$$\begin{aligned}\mathrm{E}_0[|\mathbf{R}_{i1}^{(n;\ell)}|^2] &= \frac{n_i^{-1}}{4\sigma^4}\sum_{j=1}^{n_i} \mathrm{E}_0[[(K(\tilde{G}_{1k}(d_{ij}^0/\sigma)) - k) \\ &\qquad - (K^{(\ell)}(\tilde{G}_{1k}(d_{ij}^0/\sigma)) - \mathrm{E}[K^{(\ell)}(U)])]^2] \\ &= \frac{1}{4\sigma^4}\mathrm{Var}[K(U) - K^{(\ell)}(U)] \\ &\le \frac{1}{4\sigma^4}\mathrm{E}[(K(U) - K^{(\ell)}(U))^2] \\ &= \frac{1}{4\sigma^4}\int_0^1 [K(u) - K^{(\ell)}(u)]^2\, du\end{aligned} \tag{A.14}$$

for all $n$. Clearly, $K^{(\ell)}(u)$ converges to $K(u)$, for all $u \in (0,1)$. Also, since $|K^{(\ell)}(u)|$ is bounded by $|K(u)|$, for all $\ell \geq L$, the integrand in (A.14) is bounded (uniformly in $\ell$) by $4|K(u)|^2$, which is integrable on $(0,1)$. The Lebesgue dominated convergence theorem thus implies that $\mathrm{E}_0[|\mathbf{R}_{i1}^{(n;\ell)}|^2] = o(1)$, as $\ell \to \infty$, uniformly in $n$.

(ii) Claim (ii) is the same as (i), with $d_{ij}^n/\sigma_n$ replacing $d_{ij}^0/\sigma$. Accordingly, (ii) holds under $\mathrm{P}^{(n)}_{\boldsymbol{\vartheta}^{(n)};g_1}$. That it also holds under $\mathrm{P}^{(n)}_{\boldsymbol{\vartheta};g_1}$ follows from Lemma 3.5 in [24].

(iii) Note that $|\mathcal{L}_k(K,g_1) - \mathcal{L}_k(K^{(\ell)},g_1)|^2 = |\operatorname{Cov}[K(U) - K^{(\ell)}(U), K_{g_1}(U)]|^2 \le \mathcal{L}_k(g_1) \times \operatorname{Var}[K(U) - K^{(\ell)}(U)]$, which is $o(1)$ as $\ell \to \infty$ [see (i) above]. The result then follows from the boundedness of $(s_i^{2(n)})$. □

E.C.A.R.E.S., Institut de Recherche en Statistique
and
Département de Mathématique
Université Libre de Bruxelles
Campus de la Plaine CP 210
B-1050 Bruxelles
Belgium
E-mail: mhallin@ulb.ac.be
dpaindav@ulb.ac.be
URL: http://homepages.ulb.ac.be/~dpaindav